\documentclass{article}
\usepackage[T2A]{fontenc}
\usepackage[cp1251]{inputenc}
\usepackage[english,russian]{babel}
\usepackage[tbtags]{amsmath}
\usepackage{amsfonts,amssymb,mathrsfs,amscd}
\numberwithin{equation}{section}
\newtheorem{remark}{Remark}[section]
\newtheorem{lemma}{Lemma}[section]
\newtheorem{theorem}{Theorem}[section]
\newtheorem{definition}{Definition}[section]
\newtheorem{corollary}{Corollary}[section]

\newtheorem{conclusion}{Conclusion}

\newcommand{\arccot}{\mathop{\rm arccot}}

\renewcommand{\span}{\mathop{\rm span}}
\renewcommand{\Im}{\mathop{\rm Im}}
\renewcommand{\Re}{\mathop{\rm Re}}
\newcommand{\Ker}{\mathop{\rm Ker}}

\newcommand{\rank}{\mathop{\rm rank}}
\newcommand{\col}{\mathop{\rm col}}

\newcommand{\card}{\mathop{\rm card}}

\begin{document}
\begin{Large}
\thispagestyle{empty}
\begin{center}
{\bf Inverse scattering problem for operators with a finite dimensional non-local potential\\
\vspace{5mm}
V. A. Zolotarev}\\

B. Verkin Institute for Low Temperature Physics and Engineering
of the National Academy of Sciences of Ukraine\\
47 Nauky Ave., Kharkiv, 61103, Ukraine

Department of Higher Mathematics and Informatics, V. N. Karazin Kharkov National University \\
4 Svobody Sq, Kharkov, 61077,  Ukraine

\end{center}
\vspace{5mm}

{\small {\bf Abstract.} Scattering problem for a self-adjoint integro-differential operator, which is the sum of the operator of second derivative and of a finite-dimensional self-adjoint operator, is studied. Jost solutions are found and it is shown that the scattering function has a multiplicative structure, besides, each of the multipliers is a scattering coefficient for a pair of self-adjoint operators, one of which is a one-dimensional perturbation of the other. Solution of the inverse problem is based upon the solutions to the inverse problem for every multiplier. A technique for finding parameters of the finite-dimensional perturbation via the scattering data is described.}
\vspace{5mm}


\vspace{5mm}

\begin{center}
{\bf Introduction}
\end{center}
\vspace{5mm}

Scattering theory is one of the most important areas of mathematical physics. Method of solution of inverse scattering problem (V. A. Marchenko, I. M. Gelfand and B. M. Levitan, M. G. Krein \cite{1} -- \cite{5}) play a key role in integration of non-linear partial differential equations \cite{6}. Search for the L - A pairs for non-linear equations leads to the Sturm -- Liouville operator $L$ \cite{6}. Operators with non-local (separable \cite{3}) potentials describe behavior of particles on a crystal surface and their study for different problems is given in the papers \cite{7} -- \cite{11}.

This paper is dedicated to the scattering problem for the self-adjoint operator $L_n$ in $L^2(\mathbb{R}_+)$,
$$(L_ny)(x)=-y''(x)+\int\limits_{\mathbb{R}_+}y(t)\sum\limits_{k=1}^n\alpha_k\overline{v_k(t)}v_k(x)dt$$
with the domain
$$\mathfrak{D}(L_n)=\{y(x)\in W_2^2(\mathbb{R}_+):y(0)=0\}$$
where $n\in\mathbb{N}$; $\{\alpha_k\}_1^n$ are numbers from $\mathbb{R}$; $\{v_k(x)\}_1^n$ is a set of linear independent functions such that
$$\int\limits_{\mathbb{R}_+}(1+x^2)|v_k(x)|^2dx<\infty\quad(1\leq k\leq n).$$
This paper expands the studies of paper \cite{9} on the case of $n>1$, it studies the scattering problem for a pair of self-adjoint operators $\{L_n,L_0\}$ ($L_0$ is the operator of second derivative). Note that the existence of wave operators for this pair is evident since the difference $L_n-L_0$ is finite dimensional \cite{5, 12}.

The paper consists of five sections.

Section 1 constructs Jost solutions for the operator $L_n$ and describes properties of these solutions.

Section 2 studies scattering function of a pair of self-adjoint operators $\{L_n,L_0\}$. The description of bound states of the operator $L_n$ is given and it is shown that there is a finite number of them.

Section 3 deals with the multiplicative expansion of the scattering function $S(\lambda)$ of the pair of operators $\{L_n,L_0\}$. It is established that $S(\lambda)$ is expanded into the product of multipliers $S_k(\lambda)$, besides, $S_k(\lambda)$ is given by the scattering coefficient of the pair $\{L_k,L_{k-1}\}$ (${\displaystyle L_k=L_0+\sum\limits_{p=1}^k\alpha_p\langle.,v_p\rangle v_p}$). This fact establishes relation between the finite dimensional perturbation of the operator $L_0$ and factorization of the scattering function $S(\lambda)$.

Section 4 gives the solution to the inverse scattering problem for the case of $n=1$ (see \cite{9}). Viz., for the pair $\{L_1,L_0\}$ a method of restoration of the number $\alpha_1$ and function $v_1(x)$  via the scattering coefficient $S_1(\lambda)$ is given. Scattering data are also described.

Section 5 investigates the inverse scattering problem for $n>1$. Note that the study of scattering for the pair $\{L_2,L_0\}$ leads to an interesting and important problem of $N$-expansion for Nevanlinna functions which remains unsolved. Here, we limit ourselves to the case of $L_0$-orthogonal kernels, for which the complete solution of the inverse problem is given.

\section{Jost solutions}\label{s1}

{\bf 1.1} Consider a boundary value problem on the half-axis given by the integro-differential equation
\begin{equation}
-y''(x)+\int\limits_{\mathbb{R}_+}K_n(x,t)y(t)dt=\lambda^2y(x)\quad\left(K_n(x,t)\stackrel{\rm def}{=}\sum\limits_{k=1}^n\alpha_kv_k(x)\overline{v_k}(t)\right)\label{eq1.1}
\end{equation}
and the boundary condition
\begin{equation}
y(0)=0\label{eq1.2}
\end{equation}
where $x\in\mathbb{R}_+$; $\alpha_k\in\mathbb{R}$ ($\alpha_k\not=0$, $1\leq k\leq n$, $n\in\mathbb{Z}_+$); $\lambda\in\mathbb{C}$; $\{v_k(x)\}$ is a set of non-zero linearly independent complex-valued functions on $\mathbb{R}_+$ such that
\begin{equation}
\int\limits_{\mathbb{R}_+}(1+x^2)|v_k(x)|^2dx<\infty\quad(1\leq k\leq n).\label{eq1.3}
\end{equation}

\begin{remark}\label{r1.1}
Finite dimensional perturbation in \eqref{eq1.1} is invariant relative to transform $v_k(x)\rightarrow e^{i\varphi_k}v_k(x)$ ($\varphi_k\in\mathbb{R}$, $1\leq k\leq n$). This perturbation can be normalized in two ways: (a) upon the substitution $v_k(x)\rightarrow\sqrt{|\alpha_k|}v_k(x)$, it is natural that $\alpha_k=\pm1$ ($1\leq k\leq n$); (b) we can assume that $\|v_k(x)\|_{L^2}=1$ ($1\leq k\leq n$), due to the new notation $\alpha_k\rightarrow\alpha_k\|v_k\|_{L^2}^2$ ($1\leq k\leq n$).
\end{remark}

\begin{remark}\label{r1.2}
Relation \eqref{eq1.3}, due to the Cauchy -- Bunyakovsky inequality, implies that $v_k\in L^2(\mathbb{R}_+)\cap L^1(\mathbb{R}_+)$ since
$$\int\limits_{\mathbb{R}_+}|v_k(x)|dx=\int\limits_0^1|v_k(x)|dx+\int\limits_1^\infty\frac1x|xv_k(x)|dx\leq\left\{\int\limits_0^1|v_k(x)|^2dx\right\}^{1/2}$$
$$+\left\{\int\limits_1^\infty x^2|v_k(x)|^2dx\right\}^{1/2}\leq2\left\{\int\limits_0^\infty(1+x^2)|v_k(x)|^2dx\right\}^{1/2}<\infty.$$
\end{remark}

Relation \eqref{eq1.3} implies that, for $x\rightarrow\infty$, equation \eqref{eq1.1} becomes the elementary equation $y''(x)+\lambda^2y(x)=0.$ Therefore, it is natural to search the solution $e(\lambda,x)$ to equation \eqref{eq1.1} satisfying the boundary condition
\begin{equation}
\lim\limits_{x\rightarrow\infty}e^{-i\lambda x}e(\lambda,x)=a(\lambda),\label{eq1.4}
\end{equation}
here $a(\lambda)$ is a function of $\lambda$. It is easy to show \cite{1, 2} that $e(\lambda,x)$ is the solution to the integral equation
\begin{equation}
e(\lambda,x)=e^{i\lambda x}a(\lambda)+\int\limits_x^\infty\frac{\sin\lambda(t-x)}\lambda\int\limits_{\mathbb{R}_+}e(\lambda,\tau)\sum\limits_k\alpha_k\overline{v_k(\tau)}v_k(t)d\tau dt.\label{eq1.5}
\end{equation}
Define the functions
\begin{equation}
e_k(\lambda)\stackrel{\rm def}{=}\int\limits_{\mathbb{R}_+}e(x,\lambda)\overline{v_k(x)}dx\quad(1\leq k\leq n).\label{eq1.6}
\end{equation}
Upon multiplying \eqref{eq1.5} by $\overline{v_k(x)}$ and integrating along $\mathbb{R}_+$, we obtain the system of linear equations for $\{e_k(\lambda)\}_1^n$ \eqref{eq1.6}
\begin{equation}
e_k(\lambda)=a(\lambda)\widetilde{v}_k^*+\sum\limits_{s=1}^n\alpha_se_s(\lambda)\varphi_{s,k}(\lambda)\quad(1\leq k\leq n)\label{eq1.7}
\end{equation}
where $\widetilde{v}_k(\lambda)$ is given by
\begin{equation}
\widetilde{v}_k(\lambda)\stackrel{\rm def}{=}\int\limits_{\mathbb{R}_+}e^{-i\lambda x}v_k(x)dx\quad(1\leq k\leq n)\label{eq1.8}
\end{equation}
and
\begin{equation}
\widetilde{v}_k^*(\lambda)\stackrel{\rm def}{=}\overline{\widetilde{v}_k(\overline{\lambda})}\quad(1\leq k\leq n).\label{eq1.9}
\end{equation}
The functions $\varphi_{s,k}(\lambda)$ are
\begin{equation}
\varphi_{s,k}(\lambda)\stackrel{\rm def}{=}\int\limits_{\mathbb{R}_+}\int\limits_x^\infty\frac{\sin\lambda(t-x)}\lambda v_s(t)dt\overline{v}_k(x)dx=\int\limits_{\mathbb{R}_+}\frac{\sin\lambda y}\lambda\int\limits_{\mathbb{R}_+}v_s(x+y)\overline{v_k(y)}dxdy\label{eq1.10}
\end{equation}
($1\leq s$, $k\leq n$) and $\varphi_{s,k}(-\lambda)=\varphi_{s,k}(\lambda)$. Let
\begin{equation}
g_{s,k}(y)\stackrel{\rm def}{=}\int\limits_{\mathbb{R}_+}\overline{v_s(x+y)}v_k(x)dx\quad(y\in\mathbb{R}_+,1\leq s,k\leq n).\label{eq1.11}
\end{equation}

\begin{remark}\label{r1.3}
Extend $v_k(x)$ by zero onto $\mathbb{R}_-$ ($1\leq k\leq n$), then
$$g_{s,k}(y)=\int\limits_{-y}^\infty\overline{v_s(x+y)}v_k(x)dx=\int\limits_{\mathbb{R}_+}\overline{v_s(\xi)}v_k(\xi-y)d\xi,$$
thus, $\overline{g_{s,k}(y)}=g_{k,s}(-y)$ ($1\leq s,k\leq n$). So, if $\{v_k(x)\}_1^n$ are extended by zero onto the half-axis $\mathbb{R}_-$, then the last relation gives the rule of extension of the functions $g_{s,k}(y)$ \eqref{eq1.11} onto $\mathbb{R}_-$.
\end{remark}

\begin{remark}\label{r1.4}
Convolution of two functions from $L^1(\mathbb{R})$ is a function from $L^1(\mathbb{R})$, and convolution of functions from $L^1(\mathbb{R})$ and $L^2(\mathbb{R})$ is a function from $L^2(\mathbb{R})$ \cite{13, 14}, therefore, $g_{s,k}(y)$ \eqref{eq1.11} belongs to $L^1(\mathbb{R})\cap L^2(\mathbb{R})$ because $v_k(x)\in L^1(\mathbb{R})\cap L^2(\mathbb{R})$.
\end{remark}

Equation \eqref{eq1.10} yields that
\begin{equation}
\varphi_{s,k}(\lambda)=\frac1{2i\lambda}\{\Phi_{s,k}^*(\lambda)-\Phi_{s,k}^*(-\lambda)\}\quad(1\leq s,k\leq n),\label{eq1.12}
\end{equation}
here
\begin{equation}
\Phi_{s,k}(\lambda)\stackrel{\rm def}{=}\int\limits_{\mathbb{R}_+}e^{-i\lambda y}g_{s,k}(y)dy\quad(1\leq s,k\leq n).\label{eq1.13}
\end{equation}
Since
$$\Phi_{s,k}(\lambda)=\int\limits_{\mathbb{R}_+}e^{-i\lambda y}\int\limits_{\mathbb{R}_+}\overline{v_s(x+y)}v_k(x)dxdy=\int\limits_{\mathbb{R}_+}dy\int\limits_{\mathbb{R}_+}dxe^{-i\lambda(x+y)}\overline{v_s}(x+y)e^{i\lambda x}v_k(x)$$
$$=\int\limits_{\mathbb{R}_+}e^{i\lambda x}v_k(x)\int\limits_x^\infty e^{-i\lambda t}\overline{v_s(t)}dtdx,$$
then, upon integration by parts,
\begin{equation}
\begin{array}{ccc}
{\displaystyle\Phi_{s,k}(\lambda)=-\int\limits_{\mathbb{R}_+}\int\limits_x^\infty e^{-i\lambda t}\overline{v_s(t)}dtd\left(\int\limits_x^\infty e^{i\lambda y}v_k(y)dy\right)=\int\limits_{\mathbb{R}_+}e^{-i\lambda t}\overline{v_s(t)}dt\int\limits_{\mathbb{R}_+}e^{i\lambda y}v_k(y)dy}\\
{\displaystyle-\int\limits_{\mathbb{R}_+}\int\limits_x^\infty e^{i\lambda y}v_k(y)e^{-i\lambda x}\overline{v_s}(x)dx=\widetilde{v}_s^*(-\lambda)\widetilde{v}_k(-\lambda)-\Phi_{k,s}^*(\lambda).}
\end{array}\label{eq1.14}
\end{equation}

\begin{lemma}\label{l1.1}
For the functions $\widetilde{v}_k(\lambda)$ \eqref{eq1.8} and $\Phi_{s,k}(\lambda)$ \eqref{eq1.13} ($g_{s,k}(\lambda)$ is given by \eqref{eq1.11}), for all $\lambda\in\mathbb{R}$, the following equalities are true:
\begin{equation}
\Phi_{s,k}(\lambda)+\Phi_{k,s}^*(\lambda)=\widetilde{v}_s^*(-\lambda)\widetilde{v}_k(-\lambda)\quad(1\leq s,k\leq n).\label{eq1.15}
\end{equation}
\end{lemma}

\begin{remark}\label{r1.5}
Functions $\widetilde{v_k}(\lambda)$ \eqref{eq1.8} and $\Phi_{s,k}(\lambda)$ \eqref{eq1.13} are of Hardy class $H_-^2$ \cite{13} -- \cite{15} and bounded in the closed lower half-plane $\overline{\mathbb{C}_-}$ (for $\lambda\in\mathbb{R}$, the functions $\widetilde{v}_k(\lambda)$ and $\Phi_{s,k}(\lambda)$ are bounded as Fourier transforms of functions from $L^1(\mathbb{R}_+)$ \cite{13} -- \cite{15}). Moreover, $\widetilde{v}_k(\lambda)$ and $\Phi_{s,k}(\lambda)$ are uniformly continuous when $\lambda\in\mathbb{R}$ (as Fourier transforms of functions from $L^1(\mathbb{R}_+)$) and differentiable almost everywhere, besides, $v'_k(\lambda)$, $\Phi'_{s,k}(\lambda)\in L^2(\mathbb{R})$ since $xv_k(x)\in L^2(\mathbb{R}_+)$ and $yg_{s,k}(y)\in L^2(\mathbb{R}_+)$ \cite{13, 14}. The second inclusion follows from the equality
$$yg_{s,k}(y)=\int\limits_{\mathbb{R}_+}v_k(x)\overline{v_s(x+y)}(x+y)dx-\int\limits_{\mathbb{R}_+}xv_k(x)\overline{v_s(x+y)}dx\quad(1\leq s,k\leq n)$$
and Remark \ref{r1.4}.
\end{remark}

The following statement is inverse to Lemma \ref{l1.1}.

\begin{lemma}\label{l1.2}
Let a set of functions $\{\Phi_{s,k}(\lambda)\}$ from $H_-^2$ satisfy relation \eqref{eq1.15} where $\{\widetilde{v}_k(\lambda)\}$ are Fourier transforms \eqref{eq1.8} of the functions $v_k\in L^1(\mathbb{R}_+)\cap L^2(\mathbb{R}_+)$ ($1\leq k\leq n$), then $\Phi_{s,k}(\lambda)$ are given by \eqref{eq1.13} where $\{g_{s,k}(x)\}$ are expressed via the functions $\{v_k(x)\}$ by the formulas \eqref{eq1.11}.
\end{lemma}

P r o o f. The Paley -- Wiener theorem \cite{13} -- \cite{15} implies that
$$\Phi_{s,k}(\lambda)=\int\limits_{\mathbb{R}+}e^{-i\lambda y}f_{s,k}(y)dy\quad(1\leq k,s\leq n)$$
where $f_{s,k}\in L^2(\mathbb{R}_+)$. From the functions $v_k(x)\in L^1(\mathbb{R}_+)\cap L^2(\mathbb{R}_+)$, we construct $g_{s,k}(y)$ \eqref{eq1.11} and set
$$\psi_{s,k}(\lambda)=\int\limits_{\mathbb{R}_+}e^{-i\lambda y}g_{s,k}(y)dy\quad(1\leq k,s\leq n).$$
Since for $\Phi_{s,k}$ and $\psi_{s,k}$, \eqref{eq1.15} holds, then $F_{s,k}(\lambda)+F_{k,s}^*(\lambda)=0$ where $F_{s,k}(\lambda)=\Phi_{s,k}(\lambda)-\psi_{s,k}(\lambda)$ ($1\leq s$, $k\leq n$), therefore,
$$0=\int\limits_{\mathbb{R}_+}e^{-i\lambda x}(f_{s,k}(y)-g_{s,k}(y))dy+\int\limits_{\mathbb{R}_+}e^{i\lambda x}(\overline{f_{k,s}}(y)-\overline{g_{k,s}}(y))dy=\int\limits_{\mathbb{R}_+}e^{-i\lambda y}h_{s,k}(y)dy,$$
here $h_{s,k}(y)\stackrel{\rm def}{=}(f_{s,k}(y)-g_{s,k}(y))\chi_{\mathbb{R}_+}+(\overline{f_{k,s}(y)}-\overline{g_{k,s}(y)})\chi_{\mathbb{R}_-}$ ($\chi_{\mathbb{R}_\pm}$ are characteristic functions of the sets $\mathbb{R}_\pm$). The Parseval equality implies that $h_{s,k}(y)=0$, and thus $f_{s,k}(y)=g_{s,k}(y)$ ($1\leq s$, $k\leq n$). $\blacksquare$

Study system \eqref{eq1.7}, and let
\begin{equation}
e_0(\lambda)\stackrel{\rm def}{=}\det\left[
\begin{array}{ccc}
\alpha_1\varphi_{1,1}(\lambda)-1&...&\alpha_n\varphi_{n,1}(\lambda)\\
...&...&...\\
\alpha_1\varphi_{1,n}(\lambda)&...&\alpha_n\varphi_{n,n}(\lambda)-1
\end{array}\right]\label{eq1.16}
\end{equation}
be the main determinant of system \eqref{eq1.7}. Assuming that $e_0(\lambda)\not=0$ and $a(\lambda)=e_0(\lambda)$, we obtain the solution $\{e_k(\lambda)\}_1^n$ to this system:
$$e_1(\lambda)=-\det\left[
\begin{array}{cccc}
\widetilde{v}_1^*(\lambda)&\alpha_2\varphi_{2,1}(\lambda)&...&\alpha_n\varphi_{n,1}(\lambda)\\
...&...&...&...\\
\widetilde{v}_n^*(\lambda)&\alpha_2\varphi_{2,n}(\lambda)&...&\alpha_n\varphi_{n,n}(\lambda)-1
\end{array}\right],...,$$
\begin{equation}
e_n(\lambda)=-\det\left[
\begin{array}{cccc}
\alpha_1\varphi_{1,1}(\lambda)-1&...&\alpha_{n-1}\varphi_{n-1,1}(\lambda)&\widetilde{v}_1^*(\lambda)\\
...&...&...&...\\
\alpha_1\varphi_{1,n}(\lambda)&...&\alpha_{n-1}\varphi_{n-1,n}(\lambda)&\widetilde{v}_n^*(\lambda)
\end{array}\right].\label{eq1.17}
\end{equation}

\begin{theorem}\label{t1.1}
The Jost solution $e(\lambda,x)$ to equation \eqref{eq1.1} satisfying the boundary condition \eqref{eq1.4}, where $a(\lambda)=e_0(\lambda)$ \eqref{eq1.16}, is
\begin{equation}
e(\lambda,x)=e^{i\lambda x}e_0(\lambda)+\sum\limits_{k=1}^n\alpha_ke_k(\lambda)\psi_k(\lambda,x)\label{eq1.18}
\end{equation}
where $\{e_k(\lambda)\}_1^n$ are from \eqref{eq1.17}; $e_0(\lambda)\not=0$; $\{\widetilde{v}_k(\lambda)\}$ and $\{\varphi_{s,k}(\lambda)\}$ are given by the formulas \eqref{eq1.8} and \eqref{eq1.12}; the functions $\psi_k(\lambda,x)$ are given by
\begin{equation}
\psi_k(\lambda,x)\stackrel{\rm def}{=}\int\limits_1^\infty\frac{\sin\lambda(t-x)}\lambda v_k(t)dt\quad(1\leq k\leq n).\label{eq1.19}
\end{equation}
\end{theorem}

Functions $\{\varphi_{s,k}(\lambda)\}$ and $\{\psi_k(\lambda,x)\}$ are even relative to $\lambda$, therefore the function
\begin{equation}
e(-\lambda,x)=e^{-i\lambda x}e_0(\lambda)+\sum\limits_{k=1}^n\alpha_ke_k(-\lambda)\psi_k(\lambda,x)\label{eq1.20}
\end{equation}
is the solution to equation \eqref{eq1.1} and satisfies the boundary condition
\begin{equation}
\lim\limits_{x\rightarrow\infty}e^{i\lambda x}e(-\lambda,x)=a(\lambda)(=e_0(\lambda)).\label{eq1.21}
\end{equation}

\begin{remark}\label{r1.5}
For all $p$ ($1\leq p\leq n$) and all $\lambda\in\mathbb{R}$ such that $\widetilde{v}_p(\lambda)\not=0$ and $e_0(\lambda)\not=0$, the following equality is true:
\begin{equation}
\langle e(\lambda,x)-e^{i\lambda x}(\widetilde{v}_p^*(\lambda))^{-1}e_0(\lambda),v_p(x)\rangle_{L^2}=0\quad(1\leq p\leq n).\label{eq1.22}
\end{equation}
To prove \eqref{eq1.22}, it is necessary to multiply \eqref{eq1.18} by $\overline{v}_p(x)$ and integrate along $\mathbb{R}_+$, and use the relations \eqref{eq1.7}, \eqref{eq1.10} ($a(\lambda)=e_0(\lambda)$).
\end{remark}
\vspace{5mm}

{\bf 1.2} Consider the set
\begin{equation}
E_\alpha\stackrel{\rm def}{=}\{\lambda\in\mathbb{R}:e_0(\lambda)=0\}\label{eq1.23}
\end{equation}
where $e_0(\lambda)$ is from \eqref{eq1.16}.

\begin{remark}\label{r1.7}
Set $E_\alpha$ is symmetric ($\lambda\in E_\alpha\Longleftrightarrow-\lambda\in E_\alpha)$ and bounded since, due to the Riemann -- Lebesgue lemma \cite{13} -- \cite{15}, the functions $\varphi_{s,k}(\lambda)$ \eqref{eq1.10} are small in modulo when $|\lambda|\gg1$, and thus $e_0(\lambda)=(-1)^n+\varepsilon(\lambda)$ where $|\varepsilon(\lambda)|\ll1$ ($|\lambda|\gg1$). Continuity of $\Phi_{s,k}(\lambda)$ \eqref{eq1.13} implies closedness of the set $E_\alpha$.
\end{remark}

\begin{lemma}\label{l1.3}
The functions $e(\lambda,x)$ \eqref{eq1.18} and $e(-\lambda,x)$ \eqref{eq1.20} are linearly independent if $\lambda\in\mathbb{R}\setminus E_\alpha$ ($E_\alpha$ is given by \eqref{eq1.23} and $\{v_k(x)\}$ are linearly independent).
\end{lemma}

P r o o f. Lemma's statement follows from the asymptotics \eqref{eq1.4} and \eqref{eq1.21} of the solutions $e(\lambda,x)$ \eqref{eq1.18} and $e(-\lambda,x)$ \eqref{eq1.20}. Give the straightforward proof. Assuming the contrary, we suppose that for some $\lambda\in\mathbb{R}$ there are such $\mu$, $\nu\in\mathbb{C}$ ($\mu\not=0$, $\nu\not=0$) that $\mu e(\lambda,x)+\nu e(-\lambda,x)=0$ ($\forall x\in\mathbb{R}_+$), then
$$\left\{
\begin{array}{ccc}
\mu e(\lambda,x)+\nu e(-\lambda,x)=0;\\
\mu e'(\lambda,x)+\nu e'(-\lambda,x)=0.
\end{array}\right.$$
The determinant of this system is
$$W(\lambda,x)\stackrel{\rm def}{=}e(\lambda,x)e'(-\lambda,x)-e'(-\lambda,x)e(\lambda,x)=-2i\lambda e_0^2(\lambda)+e_0(\lambda)\sum\limits_k\alpha_k(e_k(-\lambda)e^{i\lambda x}$$
$$-e_k(\lambda)e^{-i\lambda x})\psi'_k(\lambda,x)-i\lambda e_0(\lambda)\sum\limits_k\alpha_k(e_k(-\lambda)e^{i\lambda x}+e_k(\lambda)e^{-i\lambda x})\psi_k(\lambda,x)$$
$$+\sum\limits_{k,s}\alpha_k\alpha_se_k(\lambda)e_s(-\lambda)\{\psi_k(\lambda,x)\psi'_s(\lambda,x)-\psi'_k(\lambda,x)\psi_s(\lambda,x)\}.$$
Using the equalities
$$\psi'_k(\lambda,x)-i\lambda\psi_k(\lambda,x)=-\int\limits_x^\infty e^{i\lambda(t-x)}v_k(t)dt;$$
$$\psi'_k(\lambda,x)+i\lambda\psi_k(\lambda,x)=-\int\limits_x^\infty e^{-i\lambda(t-x)}v_k(t)dt,$$
we obtain that
$$W(\lambda,x)=-2i\lambda e_0^2(\lambda)+e_0(\lambda)\sum\limits_k\alpha_k\left[e_k(\lambda)\int\limits_x^\infty e^{-i\lambda t}v_k(t)dt-e_k(-\lambda)\int\limits_x^\infty e^{i\lambda t}v_k(t)dt\right]$$
$$+\sum\limits_{k,s}\alpha_k\alpha_se_k(\lambda)e_s(-\lambda)\{\psi_k(\lambda,x)\psi'_s(\lambda,x)-\psi_s(\lambda,x)\psi'_k(\lambda,x)\}.$$
If $W(\lambda,x)=0$ for some $\lambda$ (and all $x\in\mathbb{R}_+$), then $W'(\lambda,x)=0$ also. And since $\psi''_k(\lambda,x)=v_k(x)-\lambda^2\psi_k(\lambda,x)$ ($1\leq k\leq n$), then $W'(\lambda,x)=0$ implies
$$0=e_0(\lambda)\sum\limits_k\alpha_kv_k(x)(e_k(-\lambda)e^{i\lambda x}-e_k(\lambda)e^{-i\lambda x})$$
$$+\sum\limits_k\alpha_ke_k(\lambda)\psi_k(\lambda,x)\sum\limits_s\alpha_se_s(-\lambda)v_s(x)-\sum\limits_s\alpha_se_s(-\lambda)\psi_s(\lambda,x)\sum\limits_k\alpha_ke_k(
\lambda)v_k(x)$$
$$=\sum\limits_k\alpha_kv_k(x)\left\{e_k(-\lambda)\left[e_0(\lambda)e^{i\lambda x}+\sum\limits_s\alpha_se_s(\lambda)\psi_s(\lambda,x)\right]\right.$$
$$\left.-e_k(\lambda)\left[e^{-i\lambda x}e_0(\lambda)+\sum\limits_s\alpha_se_s(-\lambda)\psi_s(\lambda,x)\right]\right\},$$
and taking into account \eqref{eq1.18}, \eqref{eq1.20}, we have
$$0=\sum\limits_k\alpha_kv_k(x)\{e_k(-\lambda)e_k(\lambda,x)-e_k(\lambda)e_k(-\lambda,x)\}.$$
Multiplying this equality by $\nu$ ($\not=0$) and using $\mu e(\lambda,x)+\nu e(-\lambda,x)=0$, we obtain
$$0=e(\lambda,x)\sum\limits_k\alpha_kv_k(x)[\mu e_k(\lambda)+\nu e_k(-\lambda)]$$
which, in view of linear independence of $\{v_k(x)\}$ and $e(\lambda,x)\not=0$ gives $\mu e_k(\lambda)+\nu e_k(-\lambda)=0$ ($1\leq k\leq n$), therefore
$$0=\mu e(\lambda,x)+\nu e(-\lambda,x)=a_0(\lambda)(\mu e^{i\lambda x}+\nu e^{-i\lambda x})\quad(\forall x\in\mathbb{R}_+),$$
and thus $e_0(\lambda)=0$. $\blacksquare$

So, for $\lambda\in\mathbb{R}\setminus E_\alpha$, functions $e(\lambda,x)$ \eqref{eq1.18} and $e(-\lambda,x)$ \eqref{eq1.20} form the fundamental system of solutions to equation \eqref{eq1.1}.

Recall the well-known Titchmarsh theorem \cite{14, 16}.

\begin{theorem}[Titchmarsh]\label{t1.2}
Let $F(x)\in L^2(\mathbb{R})$, then the following statements are equivalent:

(i) the function $F(x)$ is holomorphically extendable into $\mathbb{C}_+$ and is of Hardy class $H_+^2$;

(ii) the first Sokhotski formula is true:
$$\Re F(\lambda)=\frac1\pi\int\limits_{\mathbb{R}}\hspace{-4.4mm}/\frac{\Im F(x)}{x-\lambda}dx\quad(\lambda\in\mathbb{R});$$

(iii) the second Sokhotski formula is true:
$$\Im F(\lambda)=-\frac1\pi\int\limits_{\mathbb{R}}\hspace{-4.4mm}/\frac{\Re F(x)}{x-\lambda}dx\quad(\lambda\in\mathbb{R}).$$
\end{theorem}

Integrals in (ii), (iii) are understood in the sense of principal value, and the formulas from these paragraphs are the {\bf dispersive relations} \cite{6, 16}.

Taking into account \eqref{eq1.12}, we write the function $e_0(\lambda)$ \eqref{eq1.16} as
\begin{equation}
e_0(\lambda)=\frac{\alpha_1...\alpha_n}{(2i\lambda)^n}\det\{\Phi(\lambda)-\Phi(-\lambda)-2i\lambda\alpha^{-1}\}\label{eq1.24}
\end{equation}
where
\begin{equation}
\Phi(\lambda)\stackrel{\rm def}{=}\left[
\begin{array}{ccc}
\Phi_{1,1}^*(\lambda)&...&\Phi_{n,1}^*(\lambda)\\
...&...&...\\
\Phi_{1,n}^*(\lambda)&...&\Phi_{n,n}(\lambda)
\end{array}\right];\quad\alpha\stackrel{\rm def}{=}\left[
\begin{array}{ccc}
\alpha_1&...&0\\
...&...&...\\
0&...&\alpha_n
\end{array}\right].\label{eq1.25}
\end{equation}
Relations \eqref{eq1.15} in the matrix form become
\begin{equation}
\Phi(\lambda)+\Phi^+(\lambda)=V^+(-\lambda)\quad\left(V(\lambda)\stackrel{\rm def}{=}\left[
\begin{array}{ccc}
\widetilde{v_1}(\lambda)&...&\widetilde{v_n}(\lambda)\\
...&...&...\\
0&...&0
\end{array}\right]\right),\label{eq1.26}
\end{equation}
besides, $\Phi^+(\lambda)$ and $V^+(\lambda)$ are obtained from $\Phi(\lambda)$ and $V(\lambda)$ upon the application of superposition of the operations '*' and transposition. Elements of the matrix $\Phi(\lambda)$ \eqref{eq1.25} belong to $H_+^2$ (Remark \ref{r1.5}), therefore taking into account \eqref{eq1.26} and (iii) from Theorem 1.2, we have
\begin{equation}
\Phi(\lambda)=\frac12V^+(-\lambda)V(-\lambda)+\frac1{2\pi i}\int\limits_{\mathbb{R}}\hspace{-4.4mm}/\frac{dt}{t-\lambda}V^+(-t)V(-t)\quad(\lambda\in\mathbb{R}),\label{eq1.27}
\end{equation}
and thus $\Phi(\lambda)={\it\Phi}(\lambda+i0)$ where
\begin{equation}
{\it{\Phi}}(\lambda)\stackrel{\rm def}{=}\frac1{2\pi i}\int\limits_\mathbb{R}\frac{dt}{t-\lambda}V^+(-t)V(t)\quad(t\in\mathbb{C}\setminus\mathbb{R}).\label{eq1.28}
\end{equation}
Formula \eqref{eq1.27} implies that
\begin{equation}
\begin{array}{ccc}
{\displaystyle\Phi(\lambda)-\Phi(-\lambda)-2i\lambda\alpha^{-1}=\frac12(V^+(-\lambda)V(-\lambda)-V^+(\lambda)V(\lambda))}\\
{\displaystyle-2i\lambda\left\{\frac1{2\pi}\int\limits_{\mathbb{R}}\hspace{-4.4mm}/
\frac{dt}{t^2-\lambda^2}V^+(-t)V(-t)+\alpha^{-1}\right\}.}
\end{array}\label{eq1.29}
\end{equation}

Show that for $n=1$ the set $E_\alpha$ \eqref{eq1.23} is finite. Since, for $n=1$, $\alpha=\alpha_1$ and $V^+(\lambda)V(\lambda)=|\widetilde{v}_1(\lambda)|^2$ ($\lambda\in\mathbb{R}$), then $\Phi(\lambda)-\Phi(-\lambda)-2i\lambda\alpha_1^{-1}=0$ ($\lambda\in E_1$) implies that
\begin{equation}
\left\{
\begin{array}{lll}
|\widetilde{v}_1(-\lambda)|=|\overline{v}_1(\lambda)|;\\
{\displaystyle\lambda\left\{\frac1{2\pi}\int\limits_{\mathbb{R}}\hspace{-4.4mm}/\frac{dt}{t^2-\lambda^2}|\widetilde{v}_1(-t)|^2+\frac1{\alpha_1}\right\}=0.}
\end{array}\right.\label{eq1.30}
\end{equation}
If $E_\alpha$ is infinite, then Remark \ref{r1.7} yields that there exist such sequence $\{\mu_k\}_1^\infty$ from $E_\alpha$ that $\mu_k\rightarrow\mu$ ($k\rightarrow\infty$, $\mu\in E_\alpha$) and
\begin{equation}
\frac1{2\pi}\int\limits_{\mathbb{R}}\hspace{-4.4mm}/\frac{dt}{t^2-\lambda^2}|\widetilde{v}_1(-t)|^2+\frac1{\alpha_1}=0.\label{eq1.31}
\end{equation}
Subtracting these equalities for $\lambda=\mu_k$ and $\lambda=\mu$ ($\mu_k\not=\mu$), we obtain
$$\frac1{2\pi}\int\limits_\mathbb{R}\hspace{-4.4mm}/\frac{dt}{(t^2-\mu_k^2)(t^2-\mu^2)}|\widetilde{v}_1(-t)|^2=0,$$
and upon passing to the limit $\mu_k\rightarrow\mu$,
$$\frac1{2\pi}\int\limits_{\mathbb{R}}\hspace{-4.4mm}/\frac{dt}{(t^2-\mu^2)^2}|\widetilde{v}_1(-t)|^2=0,$$
which is possible only if $\widetilde{v}_1(-t)=0$ ($\forall t\in\mathbb{R}$), but this contradicts \eqref{eq1.31}.

So, for $n=1$, the set $E_\alpha$ is
\begin{equation}
E_{\alpha_1}=\{0\}\cup\{\pm\mu_k:\mu_k>0\,(1\leq k\leq q<\infty);|\widetilde{v}_1(\mu_k)|=|\widetilde{v}(-\mu_k)|\}.\label{eq1.32}
\end{equation}

\section{Scattering function}\label{s2}

{\bf 2.1} Using $e(\lambda,x)$ \eqref{eq1.18} and $e(-\lambda,x)$ \eqref{eq1.20}, construct the function
\begin{equation}
\omega(\lambda,x)\stackrel{\rm def}{=}e(-\lambda,x)-S(\lambda)e(\lambda,x)\label{eq2.1}
\end{equation}
which is a solution to equation \eqref{eq1.1}. We define function $S(\lambda)$ from the boundary condition \eqref{eq1.2} $\omega(\lambda,0)=0$.

\begin{lemma}\label{l2.1}
For all $\lambda\in\mathbb{R}\setminus E_\alpha$ ($E_\alpha$ is given by \eqref{eq1.23}), function $\omega(\lambda,x)$ \eqref{eq2.1} is the solution to the boundary value problem \eqref{eq1.1}, \eqref{eq1.2}, besides,
\begin{equation}
S(\lambda)\stackrel{\rm def}{=}\frac{e(-\lambda,0)}{e(\lambda,0)};\quad S^{-1}(\lambda)=S(-\lambda).\label{eq2.2}
\end{equation}
\end{lemma}

The function $S(\lambda)$ \eqref{eq2.2} is said to be the {\bf scattering function} since
\begin{equation}
\omega(\lambda,x)\rightarrow e_0(\lambda)(e^{-i\lambda x}-S(\lambda)e^{i\lambda x})\quad(x\rightarrow\infty).\label{eq2.3}
\end{equation}
Using \eqref{eq1.12} and \eqref{eq1.16}, \eqref{eq1.17}, we transform $e(\lambda,0)$ \eqref{eq1.18},
$$e(\lambda,0)=e_0(\lambda)+\sum\limits_k\alpha_ke_k(\lambda)\psi_k(\lambda,0)$$
$$=\frac{\alpha_1...\alpha_n}{(2i\lambda)^n}\left\{\det\left[
\begin{array}{ccc}
{\displaystyle\Phi_{1,1}^*(\lambda)-\Phi_{1,1}^*(-\lambda)-\frac{2i\lambda}{\alpha_1}}&...&\Phi_{n,1}^*(\lambda)-\Phi_{n,1}^*(-\lambda)\\
...&...&...\\
\Phi_{1,n}^*(\lambda)-\Phi_{1,n}^*(-\lambda)&...&{\displaystyle\Phi_{n,n}^*(\lambda)-\Phi_{n,n}^*(-\lambda)-\frac{2i\lambda}{\alpha_n}}
\end{array}\right]\right.$$
$$-\det\left[
\begin{array}{ccc}
\widetilde{v}_1^*(\lambda)(\widetilde{v}_1(-\lambda)-\widetilde{v}_1(\lambda))&...&\Phi_{n,1}^*(\lambda)-\Phi_{n,1}^*(-\lambda)\\
...&...&...\\
\widetilde{v}_n^*(\lambda)(\widetilde{v}_1(-\lambda)-\widetilde{v}_1(\lambda))&...&{\displaystyle\Phi_{n,n}^*(\lambda)-\Phi_{n,n}^*(-\lambda)-\frac{2i\lambda}{\alpha_n}}
\end{array}\right]$$
$$-...-\left.\det\left[
\begin{array}{ccc}
{\displaystyle\Phi_{1,1}^*(\lambda)-\Phi_{1,1}^*(-\lambda)-\frac{2i\lambda}{\alpha_1}}&...&\widetilde{v}_1^*(\lambda)(\widetilde{v}_n(-\lambda)-\widetilde{v}_n(\lambda))\\
...&...&...\\
\Phi_{1,n}^*(\lambda)-\Phi_{1,n}^*(-\lambda)&...&\widetilde{v}_n^*(\lambda)(\widetilde{v}_n(-\lambda)-\widetilde{v}_n(\lambda))
\end{array}\right]\right\}.$$
Using \eqref{eq1.15} and proportionality of the vectors $\col[\widetilde{v}_1^*(\lambda)(\widetilde{v}_k(\lambda)-\widetilde{v}_k(\lambda)),...,\widetilde{v}_n^*(\lambda)\times$ $(\widetilde{v}_k(-\lambda)-\widetilde{v}_k(\lambda))]$ ($1\leq k\leq n$), we obtain
\begin{equation}
e(\lambda,0)=\frac{\alpha_1...\alpha_n}{(2i\lambda)^n}\label{eq2.4}
\end{equation}
$$\times\det\left[
\begin{array}{ccc}
{\displaystyle
\begin{array}{ccccc}
\Phi_{1,1}^*(\lambda)+\Phi_{1,1}(-\lambda)-\\
{\displaystyle\widetilde{v}_1^*(\lambda)\widetilde{v}_1(-\lambda)-\frac{2i\lambda}{\alpha_1}}
\end{array}}&...&{\displaystyle\begin{array}{cccc}
\Phi_{n,1}^*(\lambda)+\Phi_{1,n}(-
\lambda)-\\
\widetilde{v}_1^*(\lambda)\widetilde{v}_n(-\lambda)
\end{array}}\\
...&...&...\\
{\displaystyle\begin{array}{ccc}
\Phi_{1,n}^*(\lambda)+\Phi_{n,1}(-\lambda)-\\
\widetilde{v}_n^*(\lambda)\widetilde{v}_1(-\lambda)
\end{array}}&...&{\displaystyle\begin{array}{cccc}
\Phi_{n,n}^*(\lambda)+\Phi_{n,n}(-\lambda)\\
{\displaystyle-\widetilde{v}_n^*(
\lambda)\widetilde{v}_n(-\lambda)-\frac{2i\lambda}{\alpha_n}}
\end{array}}
\end{array}\right].$$

\begin{lemma}\label{l2.2}
For all $\lambda\in\mathbb{R}$, the following representation for the scattering function $S(\lambda)$ \eqref{eq2.2} is true:
\begin{equation}
S(\lambda)=(-1)^n\frac{r(-\lambda)}{r(\lambda)}\label{eq2.5}
\end{equation}
where
\begin{equation}
r(\lambda)\stackrel{\rm def}{=}\det R(\lambda);\quad R(\lambda)\stackrel{\rm def}{=}\Phi(\lambda)+\Phi^+(-\lambda)-V^+(\lambda)V(-\lambda)-2i\lambda\alpha^{-1};\label{eq2.6}
\end{equation}
besides, $\Phi(\lambda)$, $\alpha$, $V(-\lambda)$ are given by \eqref{eq1.25}, \eqref{eq1.26} and
\begin{equation}
\overline{r(\lambda)}=r(-\lambda);\quad\overline{S}(\lambda)=S(-\lambda);\quad S(\infty)=1.\label{eq2.7}
\end{equation}
\end{lemma}

The functions $R(\lambda)$ and $r(\lambda)$ \eqref{eq2.6} are holomorphically extendable into $\mathbb{C}_+$ (Remark \ref{r1.5}), and $R(-\lambda)$, correspondingly, into $\mathbb{C}_-$. And since $r(\lambda)\not=r(-\lambda)$ and $R(\lambda)\not=R(-\lambda)$ for $\lambda\in\mathbb{R}$, then $r(-\lambda)$ and $R(-\lambda)$ are not an analytical extension of $r(\lambda)$ and $R(\lambda)$ into $\mathbb{C}_-$.

\begin{lemma}\label{l2.3}
Matrix elements $F_{s,k}(\lambda)$ of the matrix
\begin{equation}
F(\lambda)=\Phi(\lambda)+\Phi^+(-\lambda)-V^+(\lambda)V(-\lambda)\quad(\lambda\in\mathbb{R})\label{eq2.8}
\end{equation}
belong to $H_+^2$, are bounded in the closed half-plane $\overline{\mathbb{C}_+}$ and uniformly continuous when $\lambda\in\mathbb{R}$. Functions $F_{s,k}(\lambda)$ are differentiable almost everywhere on the real axis and $F'_{s,k}(\lambda)\in L^2(\mathbb{R})$. The following equality is true:
\begin{equation}
F(\lambda)+F^+(\lambda)=[V(\lambda)-V(-\lambda)]^+[V(\lambda)-V(-\lambda)]\geq0\quad(\forall\lambda\in\mathbb{R}).\label{eq2.9}
\end{equation}
\end{lemma}

P r o o f. Boundedness, unitary continuity, and differentiability almost everywhere on $\mathbb{R}$ of the matrix elements $F_{s,k}(\lambda)$ follows from the analogous properties of the functions $\widetilde{v}_k(\lambda)$ and $\Phi_{s,k}(\lambda)$ (Remark \ref{r1.5}). Relation \eqref{eq2.9} follows from \eqref{eq1.26}:
$$F(\lambda)+F^+(\lambda)=V^+(-\lambda)V(-\lambda)+V^+(\lambda)V(\lambda)-V^+(\lambda)V(-\lambda)-V^+(-\lambda)V(\lambda)$$
$$=(V^+(\lambda)-V^+(-\lambda))(V(\lambda)-V(-\lambda)).\blacksquare$$

Applying Paragraph (iii) of Theorem \ref{t1.2} to the matrix elements $F_{s,k}(\lambda)$ and taking into account \eqref{eq2.9}, we obtain
\begin{equation}
F(\lambda)=\frac12W^+(\lambda)W(\lambda)+\frac1{2\pi i}\int\limits_\mathbb{R}\hspace{-4.4mm}/\frac{W^+(t)W(t)}{t-\lambda}dt\quad(t\in\mathbb{R})\label{eq2.10}
\end{equation}
where
\begin{equation}
W(\lambda)\stackrel{\rm def}{=}V(\lambda)-V(-\lambda)\quad(\lambda\in\mathbb{R}).\label{eq2.11}
\end{equation}
So, \eqref{eq2.10} is a matrix analogue of the Sokhotski formula \cite{17, 18} giving the boundary values $\lambda+i0$ ($\lambda\in\mathbb{R}$) on $\mathbb{R}$ from $\mathbb{C}_+$ of the Cauchy type integral
\begin{equation}
\mathcal{F}(\lambda)=\frac1{2\pi i}\int\limits_{\mathbb{R}}\frac{W^+(t)W(t)}{t-\lambda}dt\quad(\lambda\in\mathbb{C}\setminus\mathbb{R}).\label{eq2.12}
\end{equation}
Each element $F_{s,k}(-\lambda)$ belongs to $H_-^2$, then \eqref{eq2.10}, upon the substitution $t\rightarrow-t$ under the integral sign, implies
\begin{equation}
-F(-\lambda)=-\frac12W^+(\lambda)W(\lambda)+\frac1{2\pi i}\int\limits_\mathbb{R}\hspace{-4.4mm}/\frac{W^+(t)W(t)}{t-\lambda}dt\quad(\lambda\in\mathbb{R}).\label{eq2.13}
\end{equation}
The right-hand side of equality \eqref{eq2.13} coincides with the Sokhotski matrix formula \cite{17, 18} for the boundary values $\lambda-i0$ ($\lambda\in\mathbb{R}$) from $\mathbb{C}_-$ of the Cauchy type integral \eqref{eq2.12}.

\begin{remark}\label{r2.1}
Usually, the Sokhotski formulas are proved \cite{17, 18} under the supposition that the density of a Cauchy type integral satisfies the Holder condition of order no larger than $1$. The work by B. V. Khvedelidze \cite{19} implies that if the density of a Cauchy type integral belongs to $L^p(\mathbb{R})$ ($p>1$), then the special integral also belongs to $L^p(\mathbb{R})$ and non-tangential values on $\mathbb{R}$ of the Cauchy type integral exist almost everywhere and for them the Sokhotski formulas are true. In our case \eqref{eq2.10}, \eqref{eq2.12}, matrix elements of $W^+(t)W(t)$ are continuous and belong to $L^2(\mathbb{R})$, therefore the formulas \eqref{eq2.10}, \eqref{eq2.12} hold for all $\lambda\in\mathbb{R}$.
\end{remark}

\begin{lemma}\label{l2.4}
The functions
\begin{equation}
R(\lambda)=F(\lambda)-2i\lambda\alpha^{-1};\quad-R(\lambda)=-F(-\lambda)-2i\lambda\alpha^{-1}\quad(\lambda\in\mathbb{R}),\label{eq2.14}
\end{equation}
where $F(\lambda)$ is given by \eqref{eq2.8}, are the boundary values on $\mathbb{R}$ from $\mathbb{C}_+$ and $\mathbb{C}_-$ of the matrix-valued function
\begin{equation}
A(\lambda)=\frac1{2\pi i}\int\limits_{\mathbb{R}}\frac{W^+(t)W(t)}{t-\lambda}dt-2i\lambda\alpha^{-1}\quad(\lambda\in\mathbb{C}\setminus\mathbb{R}).\label{eq2.15}
\end{equation}
Moreover, $A(-\lambda)=-A(\lambda)$, and for the functions $F(\lambda)$ and $-F(-\lambda)$, the representations \eqref{eq2.10} and \eqref{eq2.13} are true.
\end{lemma}

So, the boundary values $A_\pm(\lambda)$ on $\mathbb{R}$ from $\mathbb{C}_\pm$ of the function $A(\lambda)$ \eqref{eq2.15} are
\begin{equation}
A_+(\lambda)=A(\lambda+i0)=R(\lambda);\quad A_-(\lambda)=A(\lambda-i0)=-R(-\lambda)\quad(\lambda\in\mathbb{R})\label{eq2.16}
\end{equation}
where $R(\lambda)$ is from \eqref{eq2.6}.
\vspace{5mm}

{\bf 4.2} Study the complex roots of the function
\begin{equation}
a(\lambda)\stackrel{\rm def}{=}\det A(\lambda)\label{eq2.17}
\end{equation}
where $A(\lambda)$ is given by \eqref{eq2.15}.

\begin{remark}\label{r2.2}
The set of zeros of $a(\lambda)$ \eqref{eq2.17} is bounded, closed, and symmetric ($a(\lambda)=0\Longleftrightarrow a(-\lambda)=0$) and can have limit points only at $\mathbb{R}$.

Vanishing of the integral in \eqref{eq2.15}, as $\lambda\rightarrow\infty$, implies boundedness of the set of zeros. And its closedness follows from continuity of $A(\lambda)$ in $\mathbb{C}_+$ (and in $\mathbb{C}_-$). Condition $A(-\lambda)=-A(\lambda)$ provides symmetry of the set of zeros of $a(\lambda)$, and absence of limit points of this set outside $\mathbb{R}$ is a corollary of analyticity of $A(\lambda)$ in $\mathbb{C}\setminus\mathbb{R}$.
\end{remark}

By $n_-$, we denote the number of negative elements $\{\alpha_k\}_1^n$ of the matrix $\alpha$ \eqref{eq1.25},
\begin{equation}
n_-\stackrel{\rm def}{=}\card\{\alpha_k:\alpha_k<0(1\leq k\leq n)\},\label{eq2.18}
\end{equation}
then $n_+=n-n_-$ is the quantity of the positive numbers in $\{\alpha_k\}_1^n$.

\begin{theorem}\label{t2.1}
If $n_-=0$, then function $a(\lambda)$ \eqref{eq2.17} does not vanish when $\lambda\in\mathbb{C}_+$.

If $n_->0$ and $\{v_k(x)\}_1^n$ are linearly independent, then the equation $a(\lambda)=0$ can have only finite number of roots $\{\lambda_p\}_1^m$ in $\mathbb{C}_+$ ($a(\lambda_p)=0$, $1\leq p\leq m$; $m\in\mathbb{N}$) and they all are lying on the imaginary axis $\lambda_p=i\varkappa_p$ ($\varkappa_p>0$; $1\leq p\leq m$). The subspaces
\begin{equation}
L_p\stackrel{\rm def}{=}\Ker A(\lambda_p)\quad(1\leq p\leq m)\label{eq2.19}
\end{equation}
corresponding to zeros $\lambda_p$ are such that $L_p\cap L_q=\{0\}$ for $\lambda_p\not=\lambda_q$, besides,
\begin{equation}
\sum\limits_pe_p\leq n_-\label{eq2.20}
\end{equation}
where $e_p\stackrel{\rm def}{=}\dim L_p$ ($1\leq p\leq m$).
\end{theorem}

P r o o f. Let $\lambda$ from $\mathbb{C}_+$ be a root of the equation $a(\lambda)=0$, then there exists a vector $f(\lambda)=\col[f_1(\lambda),...,f_n(\lambda)]$ (${\displaystyle\|f(\lambda)\|^2=\sum\limits_k|f_k(\lambda)|^2=1}$) such that $A(\lambda)f(\lambda)=0$. Due to \eqref{eq2.11}, \eqref{eq2.15}, it means,  that
\begin{equation}
\int\limits_\mathbb{R}\frac{\overline{w_k(t)}F(t,\lambda)}{t-\lambda}dt+\frac{4\pi}{\alpha_k}=0\quad(1\leq k\leq n),\label{eq2.21}
\end{equation}
here
\begin{equation}
F(t,\lambda)\stackrel{\rm def}{=}\sum\limits_pw_p(t)f_p(\lambda);\quad w_k(t)\stackrel{\rm def}{=}\widetilde{v}_k(t)-\widetilde{v}_k(-t)\quad(1\leq k\leq n).\label{eq2.22}
\end{equation}
Upon multiplying \eqref{eq2.21} by $\overline{f_k(\lambda)}$ and adding by $k$, we obtain
\begin{equation}
\int\limits_\mathbb{R}\frac{|F(t,\lambda)|^2}{t-\lambda}dt=-4\pi\lambda\langle\alpha^{-1}f(\lambda),f(\lambda).\label{eq2.23}
\end{equation}
If $n_-=0$, i.e., $\alpha>0$, then for $\lambda\in\mathbb{C}_+$ the left-hand side of equality \eqref{eq2.23} belongs to $\mathbb{C}_+$, and the right-hand side, to $\mathbb{C}_-$, therefore, for $n_-=0$, equation \eqref{eq2.23} does not have roots in $\mathbb{C}_+$.

Now, let $n_->0$. Since
$$\int\limits_{\mathbb{R}}\frac{|F(t,\lambda)|^2}{t-\lambda}dt=2\lambda\int\limits_{\mathbb{R}_+}\frac{|F(t,\lambda)|^2}{t^2-\lambda^2}dt,$$
due to evenness with respect to $t$ of the function $F(t,\lambda)$, then equality \eqref{eq2.23} becomes
$$\int\limits_{\mathbb{R}_+}\frac{|F(t,\lambda)|^2}{t^2-\lambda^2}dt=-2\pi\langle\alpha^{-1}f(\lambda),f(\lambda)\rangle$$
($\lambda\not=0$, in view of $\lambda\in\mathbb{C}_+$). Subtracting from this equality its complex conjugate ($\alpha=\alpha^*$), we obtain
$$(\lambda^2-\overline{\lambda}^2)\int\limits_{\mathbb{R}_+}\frac{|F(t,\lambda)|^2}{|t^2-\lambda|^2}dt=0.$$
Hence it follows that either $\lambda^2-\overline{\lambda}^2=(\lambda-\overline{\lambda})(\lambda+\overline{\lambda})=0$, and thus $\lambda$ is a purely imaginary number (because $\lambda\in\mathbb{C}_+$) or the integral vanishes which is possible only when $F(t,\lambda)=0$ ($\forall t\in\mathbb{R}_+$), and this, due to \eqref{eq2.22}, gives linear dependence of $\{w_k(t)\}$ which contradicts the linear independence of $\{v_k(x)\}$.

Let $P_\lambda$ be an orthogonal projection onto the subspace $L_\lambda=\Ker A(\lambda)$ (besides, $a(\lambda)=0$), then \eqref{eq2.15} implies
$$\int\limits_\mathbb{R}\frac{P_\lambda W^+(t)W(t)P_\lambda}{t-\lambda}+4\pi\lambda P_\lambda\alpha^{-1}P_\lambda=0,$$
and thus
$$\int\limits_{\mathbb{R}_+}\frac{P_\lambda W^+(t)W(t)P_\lambda}{t^2-\lambda^2}dt=-2\pi P_\lambda\alpha^{-1}P_\lambda.$$
The left-hand side of this equality is positive since $t^2-\lambda^2=t^2+|\lambda|^2>0$, due to the purely imaginary value of $\lambda$, and thus $P_\lambda\alpha^{-1}P_\lambda<0$. Obviously, for $l_\lambda=\rank P_\lambda$, the inequality $l_\lambda\leq n_-$ holds since the number of negative eigenvalues of the matrix $\alpha^{-1}$ equals $n_-$. So, $l_p\leq n_-$ ($l_p=\dim L_p$) for all $p$.

Let $w\in\mathbb{C}_+$ be another ($w\not=\lambda$) root of the equation $a(\lambda)=0$, and $f(w)=\col[f_1(w),...,f_n(w)]$ be the corresponding normalized by identity vector such that
\begin{equation}
\int\limits_\mathbb{R}\frac{w_k(t)F(t,w)}{t-w}dt+\frac{4\pi w}{\alpha_k}f_k(w)=0\quad(1\leq k\leq n),\label{eq2.24}
\end{equation}
here $w_k(t)$ and $F(t,w)$ are from \eqref{eq2.22}. Multiplying \eqref{eq2.21} by $\overline{f_k(w)}$ and \eqref{eq2.24}, by $\overline{f_k(\lambda)}$, and summing by $k$, we obtain
$$\int\limits_{\mathbb{R}_+}\frac{\overline{F(t,w)}F(t,\lambda)}{t^2-\lambda^2}dt+2\pi\langle\alpha^{-1}f(\lambda),f(w)\rangle=0;$$
$$\int\limits_{\mathbb{R}_+}\frac{\overline{F(t,
\lambda)}F(t,w)}{t^2-w^2}+2\pi\langle\alpha^{-1}f(w),f(\lambda)\rangle=0.$$
Subtracting from the first relation the complex conjugate to the second one, we have
$$(\lambda^2-\overline{w}^2)\int\limits_{\mathbb{R}_+}\frac{\overline{F(t,w)}F(t,\lambda)}{(t^2-\overline{w}^2)(t^2-\lambda^2)}dt=0,$$
and since $\lambda\not=\overline{w}$ ($\lambda$, $w\in\mathbb{C}_+$), then
\begin{equation}
\int\limits_{\mathbb{R}_+}\frac{\overline{F(t,w)}F(t,\lambda)}{(t^2-\overline{w}^2)(t^2-\lambda^2)}dt=0.\label{eq2.25}
\end{equation}
Hence it follows that $L_p\cap L_q=\{0\}$ for $\lambda_p\not=\lambda_q$. Really, if this is not the case and $f(\lambda_p)=f(\lambda_q)\in L_p\cap L_q$, then the last equality implies
$$\int\limits_{\mathbb{R}_+}\frac{|F(t,\lambda_p)|^2}{(t^2+\varkappa_p^2)(t^2+\varkappa_q^2)}dt=0$$
since $\lambda_p=i\varkappa_p$, $\lambda_q=i\varkappa_q$. Thus, $F(t,\lambda_q)=0$ for all $t\in\mathbb{R}_+$ which is impossible due to linear independence of the functions $\{v_k(x)\}$.

Show that the set of zeros of the function $a(\lambda)$ \eqref{eq2.17} from $\mathbb{C}_+$ is finite. Remark \ref{r2.2} and the fact that zeros of $a(\lambda)$ lie on the imaginary axis yield that zero is the only possible limit point of this set. If the set $\{\lambda_p\}$ is infinite, then $\lambda_p\rightarrow0$, besides, $\lambda=0$ is a root of $a(\lambda)$ since $A(0)=0$. Substitute $w=\lambda_p=i\varkappa_p$, $\lambda=\lambda_q=i\varkappa_q$ into \eqref{eq2.25}, then
\begin{equation}
\int\limits_{\mathbb{R}_+}\frac{\overline{F(t,\lambda_p)}F(t,\lambda_q)}{(t^2+\varkappa_p^2)(t^2+\varkappa_q^2)}dt=0.\label{eq2.26}\
\end{equation}
Vectors $f(\lambda_p)$ belong to the unit sphere in $E^\lambda$ (normalization), therefore, taking into account its compactness, we chose from the sequence $f(\lambda_p)$ a converging subsequence $f(\lambda_{p_k})$ such that $f(\lambda_{p_k})\rightarrow f(0)$ as $\lambda_{p_k}\rightarrow0$, besides, $f(0)\not=0$ ($\|f(0)\|=1$). Passing to the limit in equality \eqref{eq2.26}, as $\lambda_{p_k}\rightarrow0$, we obtain
$$\int\limits_{\mathbb{R}_+}\frac{|F(t,0)|^2}{t^4}dt=0,$$
which gives $F(t,0)=0$ for all $t\in\mathbb{R}$, and this contradicts the linear independence of $\{v_k(x)\}$.

Finally, prove that inequality \eqref{eq2.20} is true. If $\lambda_p=i\varkappa_p$ is a zero of the function $a(\lambda)$, then \eqref{eq2.21} implies
\begin{equation}
\int\limits_{\mathbb{R}_+}\frac{W^+(t)W(t)}{t^2+\varkappa_p^2}dtf(\lambda_p)=-2\pi\alpha^{-1}f(\lambda_p)\label{eq2.27}
\end{equation}
where $f(\lambda_p)\in L_p$ \eqref{eq2.19}. Consider the space
$$L=\span\{L_p:1\leq p\leq m\},$$
and let
$$f=\sum\limits_{p=1}^mf(\lambda_p)\quad(f\in L),$$
then \eqref{eq2.27} yields
$$-2\pi\alpha^{-1}f=\sum\limits_p\int\limits_{\mathbb{R}_+}\frac{W^+(t)W(t)}{t^2+\varkappa_p^2}dtf(\lambda_p).$$
Scalar multiplying this equality by $f$ and using \eqref{eq2.22}, we have
$$-2\pi\langle\alpha^{-1}f,f\rangle=\sum\limits_{p,q}\int\limits_{\mathbb{R}_+}\frac{\overline{F(t,\lambda_q)}F(t,\lambda_p)}{t^2+\varkappa_p^2}dt$$
$$=\sum\limits_{p,q}\int\limits_{\mathbb{R}_+}t^2\frac{\overline{F(t,\lambda_q)}F(t,\lambda_p)}{(t^2+\varkappa_q^2)9t^2+\varkappa_p^2)}dt+\sum\limits_{p,q}\int\limits_{\mathbb{R}_+}
\varkappa_q^2\frac{\overline{F(t,\lambda_q)}F(t,\lambda_p)}{(t^2+\varkappa_q^2)(t^2+\varkappa_p^2)}dt,$$
and, due to \eqref{eq2.26},
$$-2\pi\langle\alpha^{-1}f,f\rangle=\int\limits_{\mathbb{R}_+}t^2\sum\limits_q\frac{\overline{F(t,\lambda_q)}}{t^2+\varkappa_q^2}\sum\limits_p\frac{F(t,\lambda_p)}{t^2+\varkappa_p^2}dt+
\sum\limits_q\varkappa_q^2\int\limits_{\mathbb{R}_+}\frac{|F(t,\lambda_q)|^2}{(t^2+\varkappa_q^2)^2}dt.$$
The right-hand side of this equality is positive, and thus $\langle\alpha^{-1}f,f\rangle\leq0$, for all $f\in L$, consequently, $\dim L\leq n_-$. For every $t\in\mathbb{R}_+$, the subspace generated by the functions $F(t,\lambda_p)/t^2+\varkappa_p^2$ has the dimension $e_p$ when $f(\lambda_p)$ runs through $L_p$ \eqref{eq2.19}, in view of the linear independence of the set $\{v_k(x)\}$. Taking into account the orthogonality $F(t,\lambda_p)/t^2+\varkappa_p^2\perp F(t,\lambda_q)/t^2+\varkappa_q^2$ which takes place in $L^2(\mathbb{R}_+)$ in view of \eqref{eq2.26} when $\lambda_p\not=\lambda_q$, we obtain that to the second term of the right-hand side of the last equality there corresponds a quadratic form, matrix of which has rank
${\displaystyle\sum\limits_pe_p}$. Hence relation \eqref{eq2.21} follows. $\blacksquare$

\begin{remark}\label{r2.3}
Relation \eqref{eq2.4} implies that the function $r(\lambda)$ \eqref{eq2.6} equals
$$r(\lambda)=\frac{(2i\lambda)^n}{\alpha_1...\alpha_n}e(\lambda,0),$$
and thus $(2i\lambda)^ne(\lambda,0)$ is holomorphically extendable into $\mathbb{C}_+$, moreover, when $n_->0$ ($n_-$ from \eqref{eq2.18}), $e(i\varkappa_p,0)=0$ ($0\leq p\leq m$) where $\lambda_p=i\varkappa_p$ ($\varkappa_p>0$) are roots of the function $a(\lambda)$ \eqref{eq2.17}.
\end{remark}
\vspace{5mm}

{\bf 2.3} Taking into account \eqref{eq1.12} and \eqref{eq1.15} -- \eqref{eq1.17}, write the function $e(\lambda,x)$ \eqref{eq1.18} as
$$e(\lambda,x)=\frac{\alpha_1...\alpha_n}{(2i\lambda)^n}\left\{e^{i\lambda x}\right.$$
$$\times\det\left[
\begin{array}{ccc}
{\displaystyle\Phi_{1,1}^*(\lambda)-\Phi_{1,1}^*(-\lambda)-\frac{2i\lambda}{\alpha_1}}&...&\Phi_{n,1}^*(\lambda)-\Phi_{n,1}^*(-\lambda)\\
...&...&...\\
\Phi_{1,n}^*(\lambda)-\Phi_{1,n}^*(-\lambda)&...&{\displaystyle\Phi_{n,n}^*(\lambda)-\Phi_{n,n}^*(-\lambda)-\frac{2i\lambda}{\alpha_n}}
\end{array}\right]$$
$$-2i\lambda\psi_1(\lambda,x)\det\left[
\begin{array}{cccc}
\widetilde{v}_1^*(\lambda)&\Phi_{2,1}^*(\lambda)-\Phi_{2,1}^*(-\lambda)&...&\Phi_{n,1}^*(\lambda)-\Phi_{n,1}^*(-\lambda)\\
...&...&...&...\\
\widetilde{v}_n^*(\lambda)&\Phi_{2,n}^*(\lambda)-\Phi_{2,n}^*(-\lambda)&...&{\displaystyle\Phi_{n,n}^*(\lambda)-\Phi_{n,n}^*(-\lambda)-\frac{2i\lambda}{\alpha_n}}
\end{array}\right]-...$$
$$\left.-2i\lambda\psi_n(\lambda,x)\det\left[
\begin{array}{ccc}
{\displaystyle\Phi_{1,1}^*(\lambda)-\Phi_{1,1}^*(-\lambda)-\frac{2i\lambda}{\alpha_1}}&...&\widetilde{v}_1^*(\lambda)\\
...&...&...\\
\Phi_{1,n}^*(\lambda)-\Phi_{1,n}^*(-\lambda)&...&\widetilde{v}_n^*(\lambda)
\end{array}\right]\right\}.$$
Using the formula
$$e^{i\lambda x}\Phi_{s,k}^*(\lambda)-e^{i\lambda x}\Phi_{s,k}(-\lambda)-\widetilde{v}_k^*(\lambda)\left[\int\limits_x^\infty e^{i\lambda(t-x)}v_s(t)dt-\int\limits_x^\infty e^{-i\lambda(t-x)}v_s(t)dt\right]=b_{s,k}(\lambda,x)$$
where
\begin{equation}
b_{s,k}(\lambda,x)\stackrel{\rm def}{=}e^{i\lambda x}\Phi_{s,k}^*(\lambda)+e^{i\lambda x}\Phi_{k,s}(-\lambda)-\widetilde{v}_k^*(\lambda)\int\limits_{\mathbb{R}_+}e^{i\lambda|t-x|}v_s(t)dt-\delta_{k,s}\frac{e^{i\lambda x}2i\lambda}{\alpha_k}\label{eq2.28}
\end{equation}
($1\leq k$, $s\leq n$) and proportionality of the vectors $\col[\psi_k(\lambda,x)\widetilde{v}_1^*(\lambda),...,\psi_k(\lambda,x)\widetilde{v}_n^*(\lambda)]$ ($1\leq k\leq3$), we have
\begin{equation}
e(\lambda,x)=\frac{\alpha_1...\alpha_n}{(2i\lambda)^n}\det\left[
\begin{array}{cccc}
b_{1,1}(\lambda,x)&e^{-i\lambda x}b_{2,1}(\lambda,x)&...&e^{-i\lambda x}b_{n,1}(\lambda,x)\\
...&...&...&...\\
b_{1,n}(\lambda,x)&e^{-i\lambda x}b_{2,n}(\lambda,x)&...&e^{-i\lambda x}b_{n,n}(\lambda,x)
\end{array}\right].\label{eq2.29}
\end{equation}
Define the matrix function
\begin{equation}
B(\lambda,x)\stackrel{\rm def}{=}\left[
\begin{array}{ccc}
b_{1,1}(\lambda,x)&...&b_{n,1}(\lambda,x)\\
...&...&...\\
b_{1,n}(\lambda,x)&...&b_{n,n}(\lambda,x)
\end{array}\right]\label{eq2.30}
\end{equation}
where $\{b_{k,s}(\lambda,x)\}$ are given by \eqref{eq2.28}.

\begin{lemma}\label{l2.5}
The Jost solution $e(\lambda,x)$ \eqref{eq1.18} is expressed via the matrix $B(\lambda,x)$ \eqref{eq2.30} by the formula
\begin{equation}
e(\lambda,x)=\frac{\alpha_1...\alpha_n}{(2i\lambda)^n}e^{-i(n-1)\lambda x}\det B(\lambda,x).\label{eq2.31}
\end{equation}
\end{lemma}

Each function $b_{s,k}(\lambda,x)$ \eqref{eq2.28} ($1\leq s$, $k\leq n$), for $\lambda=i\varkappa$ ($\varkappa>0$), belongs to $L^2(\mathbb{R}_+)$ and is bounded (by $x$) since every term in \eqref{eq2.28} has this property. Using this property of $\{b_{s,k}(\lambda,x)\}$, we show that the function $e(\lambda,x)$ \eqref{eq2.31} belongs to $L^2(\mathbb{R}_+)$ and is bounded (by $x$) when $\lambda=i\varkappa$ ($\varkappa>0$). Rewrite $b_{s,k}(\lambda,x)$ as
\begin{equation}
b_{s,k}(\lambda,x)=e^{i\lambda x}a_{s,k}(\lambda)-w_s(\lambda,x)\widetilde{v}_k^*(\lambda)\quad(1\leq s,k\leq n)\label{eq2.32}
\end{equation}
where
\begin{equation}
a_{s,k}(\lambda)\stackrel{\rm def}{=}\Phi_{s,k}^*(\lambda)+\Phi_{k,s}(-\lambda)-\delta_{s,k}\frac{2i\lambda}{\alpha_k};\quad w_s(\lambda,x)\stackrel{\rm def}{=}\int\limits_{\mathbb{R}_+}e^{i\lambda|t-x|}v_s(t)dt\label{eq2.33}
\end{equation}
($1\leq s$, $k\leq n$), then \eqref{eq2.29} implies
$$e(\lambda,x)=\frac{\alpha_1...\alpha_n}{(2i\lambda)^n}\left\{\det\left[
\begin{array}{cccc}
b_{1,1}(\lambda,x)&a_{2,1}(\lambda)&...&e^{-i\lambda x}b_{n,1}(\lambda,x)\\
...&...&...&...\\
b_{1,n}(\lambda,x)&a_{2,n}(\lambda)&...&e^{-i\lambda x}b_{n,n}(\lambda,x)
\end{array}\right]\right.$$
$$\left.-e^{-i\lambda x}w_2(\lambda,x)\det\left[
\begin{array}{cccc}
b_{1,1}(\lambda,x)&\widetilde{v}_1^*(\lambda)&...&e^{-i\lambda x}b_{n,1}(\lambda,x)\\
...&...&...&...\\
b_{1,n}(\lambda,x)&\widetilde{v}_n^*(\lambda)&...&e^{-i\lambda x}b_{n,n}(\lambda,x)
\end{array}\right]\right\}.$$
Using \eqref{eq2.32} and proportionality of the vectors $\col[w_s(\lambda,x)\widetilde{v}_1^*(\lambda),...,w_s(\lambda,x)\widetilde{v}_n^*(\lambda)]$ ($1\leq s\leq n$), we obtain that
$$e(\lambda,x)=\frac{\alpha_1...\alpha_n}{(2i\lambda)^n}\left\{\det\left[
\begin{array}{cccc}
b_{1,1}(\lambda,x)&a_{2,1}(\lambda)&...&a_{n,1}(\lambda)\\
...&...&...&...\\
b_{1,n}(\lambda,x)&a_{2,n}(\lambda)&...&a_{n,n}(\lambda)
\end{array}\right]\right.$$
$$-w_2(\lambda,x)\det\left[
\begin{array}{cccc}
a_{1,1}(\lambda)&\widetilde{v}_1^*(\lambda)&...&a_{n,1}(\lambda)\\
...&...&...&...\\
a_{1,n}(\lambda)&\widetilde{v}_n^*(\lambda)&...&a_{n,n}(\lambda)
\end{array}\right]-...$$
$$\left.-w_n(\lambda,x)\det\left[
\begin{array}{cccc}
a_{1,1}(\lambda)&...&a_{n-1,1}(\lambda)&\widetilde{v}_1^*(\lambda)\\
...&...&...&...\\
a_{1,n}(\lambda)&...&a_{n-1,n}(\lambda)&\widetilde{v}_n^*(\lambda)
\end{array}\right]\right\}.$$
Since $b_{1,k}(\lambda,x)$ and $w_k(\lambda,x)$ ($1\leq k\leq n$) belong to $L^2(\mathbb{R}_+)$ as functions of $x$ for every $\lambda=i\varkappa$ ($\varkappa>0$) and $a_{k,s}(i\varkappa)$, $\widetilde{v}_k^*(i\varkappa)$ ($1\leq k$, $s\leq n$) are bounded, then hence it follows that $e(i\varkappa,x)\in L^2(\mathbb{R}_+)$.

\begin{lemma}\label{l2.6}
Let $n_->0$ ($n_-$ is given by \eqref{eq2.18}) and $\lambda_p=i\varkappa_p$ ($\varkappa_p>0$, $1\leq p\leq m$) are zeros of the function $a(\lambda)$ \eqref{eq2.17}, then the Jost solutions $e(\lambda_p,x)$ \eqref{eq1.18}, \eqref{eq2.31} of the boundary value problem \eqref{eq1.1}, \eqref{eq1.2} belong to $L^2(\mathbb{R}_+)$.
\end{lemma}

Remark \ref{r2.3} implies validity of the boundary condition $e(\lambda_p,0)=0$.

\begin{corollary}\label{c2.1}
For $n_->0$ ($n_-$ is from \eqref{eq2.18}), there exist no more than $n_-$ {\bf bound states} \cite{1, 6, 20} which specify the eigenfunctions $e(i\varkappa_p,x)$ ($\in L^2(\mathbb{R}_+)$ and are given by \eqref{eq2.31}) of the boundary value problem \eqref{eq1.1}, \eqref{eq1.2} corresponding to the eigenvalues $\lambda_p^2=-\varkappa_p^2$ ($\lambda_p=i\varkappa_p$, $\varkappa_p>0$ are zeros of the function $a(\lambda)$ \eqref{eq2.17}).
\end{corollary}
\vspace{5mm}

{\bf 2.4} Proceed to the real zeros of the function $r(\lambda)$ \eqref{eq2.6}. Equation $r(\lambda)=0$ ($\lambda\in\mathbb{R}$) implies that there exists such a vector $f(\lambda)=\col[f_1(\lambda),...,f_n(\lambda)]$ ($\|f(\lambda)\|^2={\displaystyle\sum|f_k(\lambda)|^2=1}$) that $R(\lambda)f(\lambda)=0$ where $R(\lambda)$ is given by \eqref{eq2.6}. Since $R(\lambda)=F(\lambda)-2i\lambda^{-1}$ ($F(\lambda)$ is from \eqref{eq2.8}), then, taking into account \eqref{eq2.10}, we obtain
$$F(\lambda)=\frac12W^+(\lambda)W(\lambda)+\frac{\lambda}{\pi i}\int\limits_{\mathbb{R}_+}\hspace{-4.4mm}/\frac{W^+(t)W(t)}{t^2-\lambda^2}dt,$$
and thus
\begin{equation}
\frac12W^+(\lambda)W(\lambda)f(\lambda)-2i\lambda\left\{\frac1{2\pi}\int\limits_{\mathbb{R}_+}\hspace{-4.4mm}/\frac{W^+(t)W(t)}{t^2-\lambda^2}dt+\alpha^{-1}\right\}f(\lambda)=0.
\label{eq2.34}
\end{equation}
Scalar multiplying this equality by $f(\lambda)$ and equating to zero real and imaginary parts ($\lambda\in\mathbb{R}$) of the obtained, we have
$$f^*(\lambda)W^+(\lambda)W(\lambda)f(\lambda)=0;\quad f^*(\lambda)\left\{\int\limits_{\mathbb{R}_+}\hspace{-4.4mm}/\frac{W^+(t)W(t)}{t^2-\lambda^2}dt+2\pi\alpha^{-1}\right\}f(\lambda)=0.$$
The first equality implies $W(\lambda)f(\lambda)=0$, then using \eqref{eq2.34} we arrive at the following statement.

\begin{remark}\label{r2.4}
If $\lambda$ ($\in\mathbb{R}$) is a zero of the function $r(\lambda)$ \eqref{eq2.6}, then there exists a vector $f(\lambda)\in E^n$ ($\|f(\lambda)=1$) such that
\begin{equation}
\left\{
\begin{array}{lll}
W_(\lambda)f(\lambda)=0;\\
{\displaystyle\left(\int\limits_{\mathbb{R}_+}\hspace{-4.4mm}/\frac{W^+(t)W(t)}{t^2-\lambda^2}dt+2\pi\alpha^{-1}\right)f(\lambda)=0.}
\end{array}\right.\label{eq2.35}
\end{equation}
where $W(\lambda)$ is given by \eqref{eq2.11}.
\end{remark}

\begin{lemma}\label{l2.7}
If functions $\{v_k(x)\}$ are linearly independent, then the set of the real zeros of $r(\lambda)$ \eqref{eq2.6} is finite,
\begin{equation}
E_r\stackrel{\rm def}{=}\{0,\pm\lambda_k:r(\lambda_k)=0\,(1\leq k\leq q,\,q\in\mathbb{N})\}.\label{eq2.36}
\end{equation}
\end{lemma}

P r o o f. The set of real zeros of the function $r(\lambda)$ is bounded, closed and symmetric ($r(\lambda)=0\Longleftrightarrow r(-\lambda)=0$), due to \eqref{eq2.7}. Closedness and boundedness follow from the continuity of the matrix elements of $F(\lambda)$ \eqref{eq2.8} and the fact that $F(\lambda)\rightarrow0$ ($\lambda\rightarrow\infty$), in view of the Riemann -- Lebesgue lemma.

If this set of zeros of $r(\lambda)$ is infinite, then hence it follows that there exists a converging sequence $\lambda_k\rightarrow w$ and $r(\lambda_k)=0$ ($\forall k$), and also $r(w)=0$. Taking into account the compactness of the unit sphere in $E^r$, we select from $f(\lambda_k)$ (satisfying \eqref{eq2.34}) a converging subsequence $f(\widetilde{\lambda}_k)\rightarrow f(w)$ when $\widetilde{\lambda}_k\rightarrow w$. Using Remark \ref{r2.4}, we obtain
$$\left\{\int\limits_{\mathbb{R}_+}\hspace{-4.4mm}/\frac{W^+(t)W(t)}{t^2-\widetilde{\lambda}_k^2}dt+2\pi\alpha^{-1}\right\}f(\widetilde{\lambda}_k)=0\quad(\forall k),$$
and thus
$$f^*(\widetilde{\lambda}_q)\left\{\int\limits_{\mathbb{R}_+}\hspace{-4.4mm}/\frac{W^+(t)W(t)}{t^2-\widetilde{\lambda}_k^2}dt+2\pi\alpha^{-1}\right\}f(\widetilde{\lambda}_k)=0\quad(
\forall k,q).$$
Analogously,
$$f^*(\widetilde{\lambda_k})\left\{\int\limits_{\mathbb{R}_+}\hspace{-4.4mm}/\frac{W^+(t)W(t)}{t^2-\widetilde{\lambda}_q^2}dt+2\pi\alpha^{-1}\right\}f(\widetilde{\lambda}_q)=0\quad(
\forall q,k).$$
Subtracting from the first equality the complex conjugate of the second, we obtain
$$f^*(\widetilde{\lambda_q})\left\{\int\limits_{\mathbb{R}_+}\hspace{-4.4mm}/\frac{W^+(t)W(t)}{(t^2-\widetilde{\lambda}_k^2)(t^2-\widetilde{\lambda}_q^2)}dt\right\}f(\widetilde{
\lambda_k}=0\quad(\forall k,q),$$
since $\widetilde{\lambda}_k^2\not=\widetilde{\lambda}_q^2$. Upon passing to the limit $\widetilde{\lambda}_k$, $\widetilde{\lambda}_q\rightarrow w$, we have
$$f^*(w)\left\{\int\hspace{-4.4mm}/\frac{W^+(t)W(t)}{(t^2-w^2)^2}dt\right\}f(w)=0,$$
which is possible only under condition $W(t)f(w)=0$ ($\forall t\in\mathbb{R}_+$), and this means that the functions $\{v_k(x)\}$ are linearly dependent. $\blacksquare$

Let $\lambda_k\in E_r$ \eqref{eq2.36}, then there exists such vector $f(\lambda_k)\in E^n$ ($\|f(\lambda_k)\|=1$) that $W(\lambda_k)f(\lambda_k)=0$ and $R(\lambda_k)f(\lambda_k)=0$. Taking \eqref{eq2.6}, \eqref{eq2.11} into account, we obtain
$$0=R(\lambda_k)f(\lambda_k)=\{\Phi(\lambda_k)+\Phi^+(-\lambda_k)-V^+(\lambda_k)V(-\lambda_k)-2i\lambda_k\alpha^{-1}\}f(\lambda_k)$$
$$=\{\Phi(\lambda_k)+\Phi^+(-\lambda_k)-V^+(\lambda_k)V(\lambda_k)-2i\lambda_k\alpha^{-1}\}f(\lambda_k)$$
$$=\{\Phi(\lambda_k)-\Phi(-\lambda_k)-2i\lambda_k\alpha^{-1}\}f(\lambda_k),$$
due to \eqref{eq1.26}. Therefore, $e_0(\lambda_k)=0$ (see \eqref{eq2.4}) and solution $e(\lambda_k,x)$ \eqref{eq1.18} equals
\begin{equation}
e(\lambda_k,x)=\sum\limits_p\alpha_pe_p(\lambda_k)\psi_p(\lambda_k,x)\label{eq2.37}
\end{equation}
where $\psi_p(\lambda,x)$ is given by \eqref{eq1.19}. Function $e(\lambda_k,x)$ is the solution to the boundary value problem \eqref{eq1.1}, \eqref{eq1.2}. To show that $e(\lambda_k,x)\in L^2(\mathbb{R}_+)$, we use a theorem by G. Hardy \cite{14, 21}.

\begin{theorem}[Hardy]\label{t2.2}
If $f\in L^p(\mathbb{R}_+)$ ($p>1$), then functions
$$\varphi(x)=\frac1\lambda\int\limits_0^xf(t)dt;\quad\psi(x)=\int\limits_x^\infty\frac{f(t)}tdt$$
also belong to $L^p(\mathbb{R}_+)$.
\end{theorem}

Every function $\psi_p(\lambda_k,x)$ in \eqref{eq2.37} belongs to $L^2(\mathbb{R}_+)$ since
$$|\psi_p(\lambda_k,x)|\leq\frac1{|\lambda_k|}\int\limits_x^\infty|v_p(t)|dt=\frac1{|\lambda_k|}\int\limits_x^\infty\frac{|tv_p(t)|}tdt,$$
then, taking into account that $tv_p(t)\in L^2(\mathbb{R}_+)$ \eqref{eq1.3}, we obtain that $\psi_p(\lambda_k,x)\in L^2(\mathbb{R})$.

\begin{lemma}\label{l2.8}
For all $\lambda_k$ ($\not=0$) from $E_r$ \eqref{eq2.36}, the function $e(\lambda_k,x)$ \eqref{eq2.37} is a solution to the boundary value problem \eqref{eq1.1}, \eqref{eq1.2} and belongs to $L^2(\mathbb{R}_+)$.
\end{lemma}

\begin{corollary}\label{c2.3}
If the set $E_r$ \eqref{eq2.36} is not empty ($q>1$), then there exist a finite number of {\bf bound states} \cite{1, 6, 16} which specify the eigenfunctions $e(\lambda_k,x)$ \eqref{eq2.37} corresponding to the eigenvalues $\lambda_k^2$ ($>0$) where $\lambda_k$ are real zeros of $r(\lambda)$ \eqref{eq2.6}.
\end{corollary}

\section{Multiplicative expansion of the scattering function}\label{s3}

{\bf 3.1} Consider the self-adjoint operator $L_0$ in $L^2(\mathbb{R}_+)$,
\begin{equation}
(L_0y)(x)\stackrel{\rm def}{=}-y'(x)\label{eq3.1}
\end{equation}
with the domain
\begin{equation}
\mathfrak{D}(L_0)\stackrel{\rm def}{=}\{y(x)\in W_2^2(\mathbb{R}_+):y(0)=0\}.\label{eq3.2}
\end{equation}
Following \cite{4}, by $\varphi(\lambda,x)$ and $\theta(\lambda,x)$ we denote the functions
$$\left\{
\begin{array}{lll}
-\varphi''(\lambda,x)=\lambda^2\varphi(\lambda,x);\\
\varphi(\lambda,0)=0;\,\varphi'(\lambda,0)=1;
\end{array}\right.\quad\left\{
\begin{array}{lll}
-\theta''(\lambda,x)=\lambda^2\theta(\lambda,x);\\
\theta(\lambda,0)=1;\,\theta'(0,\lambda)=0,
\end{array}\right.$$
then $\varphi(\lambda,x)=\sin\lambda x/\lambda$; $\theta(\lambda,x)=\cos\lambda x$. For $\lambda\in\mathbb{C}\setminus\mathbb{R}$, there exist \cite{4} such function $n(\lambda)$ ($m(z)=n(\sqrt{z})$ is the Weyl function \cite{4}) that
$$\psi(\lambda,x)=\theta(\lambda,x)+n(\lambda)\varphi(\lambda,x)\in L^2(\mathbb{R}_+)$$
and $\psi(\lambda,x)$ is said to be the Weyl solution \cite{4}. For the operator $L_0$ \eqref{eq3.1}, \eqref{eq3.2}, the functions $n(\lambda)$ and $\psi(\lambda,x)$ are
$$n(\lambda)=\pm i\lambda\,(\lambda\in\mathbb{C}_\pm);\quad\psi(\lambda,x)=e^{\pm i\lambda x}\,(\lambda\in\mathbb{C}_\pm).$$
Resolvent $R(\lambda^2)$ of the Sturm -- Liouville operator is expressed \cite{4} via $\varphi(\lambda,x)$, $\psi(\lambda,x)$ by the formula
$$(R(\lambda^2)f)(x)=\psi(\lambda,x)\int\limits_0^x\varphi(\lambda,y)f(y)dy+\varphi(\lambda,x)\int\limits_x^\infty\psi(\lambda,y)f(y)dy\quad(\lambda\in\mathbb{C}\setminus\mathbb{R}),$$
therefore, for the resolvent $R_0(\lambda^2)=(L_0-\lambda^2I)^{-1}$ of operator $L_0$ \eqref{eq3.1}, \eqref{eq3.2}, we have
\begin{equation}
(R_0(\lambda^2)f)(x)=e^{\pm i\lambda x}\int\limits_0^x\frac{\sin\lambda y}\lambda f(y)dy+\frac{\sin\lambda x}\lambda\int\limits_x^\infty e^{\pm i\lambda y}f(y)dy\quad(\lambda\in\mathbb{C}_\pm).\label{eq3.3}
\end{equation}

\begin{lemma}\label{l3.1}
For all $\lambda\in\mathbb{C}_+$, the following equalities hold:
\begin{equation}
\langle R_0(\lambda^2)v_k,v_s\rangle=\frac1{2i\lambda}\{\widetilde{v}_k(-\lambda)\widetilde{v}_s^*(\lambda)-\Phi_{k,s}^*(\lambda)-\Phi_{s,k}(-\lambda)\}\label{eq3.4}
\end{equation}
($1\leq k$, $s\leq n$) where $\{\widetilde{v}_k(\lambda)\}$ and $\{\Phi_{s,k}(\lambda)\}$ are given by \eqref{eq1.8} and \eqref{eq1.13}.
\end{lemma}

P r o o f. Equation \eqref{eq3.3}, for $\lambda\in\mathbb{C}_+$ implies
$$2i\lambda\langle R_0(\lambda^2)v_k,v_s\rangle=\int\limits_0^\infty e^{i\lambda x}\overline{v_s(x)}dx\int\limits_0^xe^{i\lambda y}v_k(y)dy-\int\limits_0^\infty e^{i\lambda x}\overline{v_s}(x)dx\int\limits_0^xe^{-i\lambda y}v_k(y)dy$$
$$+\int\limits_0^\infty e^{i\lambda x}\overline{v_s(x)}dx\int\limits_x^\infty e^{i\lambda y}v_k(y)dy-\int\limits_0^\infty e^{-i\lambda x}\overline{v_s}(x)dx\int\limits_x^\infty e^{i\lambda y}v_k(y)dy$$
$$=\int\limits_0^\infty\left(\int\limits_0^xe^{i\lambda t}\overline{v_s(t)}dt\int\limits_0^xe^{i\lambda y}v_k(y)dy\right)'-\int\limits_0^\infty\int\limits_y^\infty e^{i\lambda x}\overline{v_s(x)}dxe^{-i\lambda y}v_k(y)dy$$
$$-\int\limits_0^\infty e^{i\lambda x}\overline{v_s(x)}dx\int\limits_x^\infty e^{i\lambda y}v_k(y)dy=\widetilde{v}_k(-\lambda)\widetilde{v}_s^*(\lambda)-\int\limits_0^\infty v_k(y)dy\int\limits_0^\infty e^{i\lambda\xi}\overline{v_s(\xi+y)}d\xi$$
$$-\int\limits_0^\infty\overline{v_s}(x)dx\int\limits_0^\infty e^{i\lambda\xi}v_k(x+\xi)d\xi=\widetilde{v}_k(-\lambda)\widetilde{v}_s^*(\lambda)-\Phi_{k,s}^*(\lambda)-\Phi_{s,k}(-\lambda),$$
due to \eqref{eq1.11}. $\blacksquare$

\begin{remark}\label{r3.1}
For all $\lambda\in\mathbb{C}_-$, analogously to \eqref{eq3.4},
\begin{equation}
\langle R_0(\lambda^2)v_k,v_s\rangle=-\frac1{2i\lambda}\{\widetilde{v}_k(\lambda)\widetilde{v}_s^*(-\lambda)-\Phi_{k,s}^*(-\lambda)-\Phi_{s,k}(\lambda)\}\label{eq3.5}
\end{equation}
($1\leq k$, $s\leq n$). So, when passing from $\mathbb{C}_+$ to $\mathbb{C}_-$ in equalities \eqref{eq3.4}, one has to substitute $\lambda\rightarrow-\lambda$.
\end{remark}

Define the matrix function
\begin{equation}
T(z)\stackrel{\rm def}{=}\left[
\begin{array}{cccc}
\langle R_0(z)v_1,v_1\rangle&...&\langle R_1(z)v_n,v_1\rangle\\
...&...&...\\
\langle R_1(z)v_1,v_n\rangle&...&\langle R_n(z)v_n,v_n\rangle
\end{array}\right]\label{eq3.6}
\end{equation}
where $R_0(z)=(L_0-zI)^{-1}$ and $z\in\mathbb{C}\setminus\mathbb{R}$, then, taking into account \eqref{eq3.4}, we obtain that
\begin{equation}
I+\alpha T(\lambda^2)=-\frac\alpha{2i\lambda}\{\Phi(\lambda)+\Phi^{+}(-\lambda)-V^+(\lambda)V(-\lambda)-2i\lambda\alpha^{-1}\},\label{eq3.7}
\end{equation}
due to \eqref{eq1.25}, \eqref{eq1.26}.

\begin{theorem}\label{t3.1}
For the function $e(\lambda,0)$ \eqref{eq2.4}, for all $\lambda\in\mathbb{C}_+$, the following representation holds:
\begin{equation}
e(\lambda,0)=(-1)^n\det(I+\alpha T(\lambda^2))\label{eq3.8}
\end{equation}
where $T(z)$ is from \eqref{eq3.6}, $\alpha$, from \eqref{eq1.25}.
\end{theorem}

\begin{remark}\label{r3.2}
For $e(-\lambda,0)$ when $\lambda\in\mathbb{C}_-$,
\begin{equation}
e(-\lambda,0)=(-1)^n\det(I+\alpha T(\lambda^2)),\label{eq3.9}
\end{equation}
and for the elements $\langle R_0(\lambda^2)v_k,v_s\rangle$ of the matrix $T(\lambda^2)$ \eqref{eq3.6}, representations \eqref{eq3.5} hold.

So, for $\lambda\in\mathbb{R}$, the functions $e(\lambda,0)$ (and $r(\lambda)$ \eqref{eq2.6}) and $e(-\lambda,0)$ ($r(-\lambda)$) are the boundary values on $\mathbb{R}$ from $\mathbb{C}_+$ and from $\mathbb{C}_-$ of the function $(-1)^n\det(I+\alpha T(\lambda^2))$.
\end{remark}
\vspace{5mm}

{\bf 3.2} Define the self-adjoint operators $L_k$ in $L^2(\mathbb{R}_+)$,
\begin{equation}
L_k\stackrel{\rm def}{=}L_0+\sum\limits_{s=1}^k\alpha_s\langle.,v_s\rangle v_s\quad(0\leq k\leq n)\label{eq3.10}
\end{equation}
where $L_0$ is given by \eqref{eq3.1}, \eqref{eq3.2}; $\alpha_k\in\mathbb{R}$ ($1\leq k\leq n$); $\{v_k\}_1^n$ are linearly independent functions satisfying condition \eqref{eq1.3}. By $b(z)$, we denote a scalar function,
\begin{equation}
b(z)\stackrel{\rm def}{=}\det(I+\alpha T(z))\quad(z\in\mathbb{C}\setminus\mathbb{R})\label{eq3.11}
\end{equation}
where $\alpha$ and $T(z)$ are from \eqref{eq1.25} and \eqref{eq3.6}.

\begin{theorem}\label{t3.2}
For all $z\in\mathbb{C}\setminus\mathbb{R}$, the function $b(z)$ \eqref{eq3.11} is expanded into the product
\begin{equation}
b(z)=b_1(z)...b_n(z)\label{eq3.12}
\end{equation}
where
\begin{equation}
b_k(z)\stackrel{\rm def}{=}1+\alpha_k\langle R_{k-1}(z)v_k,v_k\rangle\quad(1\leq k\leq n),\label{eq3.13}
\end{equation}
besides, $R_k(z)=(L_k-zI)^{-1}$ is resolvent of the operator $L_k$ \eqref{eq3.10} ($1\leq k\leq n$).
\end{theorem}

P r o o f. Resolvent of the operator $L_1$ \eqref{eq3.10} is \cite{10, 12}
\begin{equation}
R_1(z)f=R_0(z)f-\frac{\alpha_1\langle R_0(z)f,v_1\rangle}{1+\alpha_1\langle R_0(z)v_1,v_1\rangle}R_0(z)v_1\label{eq3.14}
\end{equation}
where $f\in L^2(\mathbb{R}_+)$. Assuming that $f=v_2$ and scalar multiplying \eqref{eq3.14} by $\alpha_2v_2$, we obtain
$$(1+\alpha_2\langle R_1(z)v_2,v_2\rangle)(1+\alpha_1\langle R_0(z)v_1,v_1\rangle)=(1+\alpha_2\langle R_0(z)v_2,v_2\rangle)(1+\alpha_1\langle R_0(z)v_1,v_1\rangle)$$
$$-\alpha_1\langle R_0(z)v_1,v_2\rangle\cdot\alpha_2\langle R_0(z)v_2,v_1\rangle=b(z)$$
where $b(z)$ coincides with \eqref{eq3.11} when $n=2$. So, \eqref{eq3.12} for $n=2$ is proved.

In view of mathematical induction, let equation \eqref{eq3.12} hold for $n$, prove its validity for $n+1$. Write the operators $L_k$ \eqref{eq3.10} ($0\leq k\leq n+1$) as
$$L_k=L_1+\sum\limits_{s=2}^k\alpha_s\langle .,v_s\rangle v_s\quad(1\leq k\leq n+1).$$
For this set of operators, due to induction assumption,
\begin{equation}
\widetilde{b}(z)=b_2(z)...b_{n+1}(z)\label{eq3.15}
\end{equation}
where
$$\widetilde{b}(z)\stackrel{\rm def}{=}\det\left[
\begin{array}{cccc}
\alpha_2\langle R_1(z)v_2,v_2\rangle+1&...&\alpha_2\langle R_1(z)v_{n+1},v_2\rangle\\
...&...&...\\
\alpha_{n+1}\langle R_1(z)v_2,v_{n+1}\rangle&...&\alpha_{n+1}\langle R_1(z)v_{n+1},v_{n+1}\rangle+1
\end{array}\right].$$
Formula \eqref{eq3.14} implies
$$\alpha_s\langle R_1(z)v_k,v_s\rangle+\delta_{k,s}=\frac{a_{s,k}a_{1,1}-a_{s,1}a_{1,k}}{a_{1,1}}\quad(2\leq s,k\leq n+1),$$
here
$$a_{s,k}=\alpha_s\langle R_0(z)v_k,v_s\rangle+\delta_{k,s}\quad(1\leq k,s\leq n+1).$$
Therefore
$$\widetilde{b}(z)=\frac1{a_{1,1}^n}\det\left[
\begin{array}{ccc}
a_{2,2}a_{1,1}-a_{2,1}a_{1,2}&...&a_{2,n+1}a_{1,1}-a_{2,1}a_{1,n+1}\\
...&...&...\\
a_{n+1,2}a_{1,1}-a_{n+1,1}a_{1,2}&...&a_{n+1,n+1}a_{1,1}-a_{n+1,1}a_{1,n+1}
\end{array}\right]$$
$$=\frac1{a_{1,1}^{n-1}}\det\left[
\begin{array}{ccccc}
a_{2,2}&a_{2,3}a_{1,1}-a_{2,1}a_{1,3}&...&a_{2,n+1}a_{n,n}-a_{2,1}a_{1,n+1}\\
...&...&...&...\\
a_{n+1,2}&a_{n+1,3}a_{1,1}-a_{n+1,1}a_{1,3}&...&a_{n+1,n+1}a_{1,1}-a_{n+1,1}a_{1,n+1}
\end{array}\right]$$
$$-\frac{a_{1,2}}{a_{2,1}}\det\left[
\begin{array}{cccc}
a_{2,1}&a_{2,3}&...&a_{2,n+1}\\
...&...&...&...\\
a_{n+1,1}&a_{n+1,3}&...&a_{n+1,n+1}
\end{array}\right],$$
due to proportionality of $\col[a_{2,1}a_{1,s},...,a_{n+1,1}a_{1,s}]$ ($2\leq s\leq n+1$). Upon repeating this procedure for other columns of the first determinant, we obtain
$$a_{1,1}\widetilde{b}_1(z)=a_{1,1}\det\left[
\begin{array}{cccc}
a_{2,2}&...&a_{2,n+1}\\
...&...&...\\
a_{n+1,2}&...&a_{n+1,n+1}
\end{array}\right]-a_{1,2}\det\left[
\begin{array}{ccc}
a_{2,1}&...&a_{2,n+1}\\
...&...&...\\
a_{n+1,1}&...&a_{n+1,n+1}
\end{array}\right]+...$$
$$+(-1)^na_{1,n+1}\det\left[
\begin{array}{ccc}
a_{2,1}&...&a_{2,n}\\
...&...&...\\
a_{n+1,1}&...&a_{n+1,n}
\end{array}\right]=\det\left[
\begin{array}{ccc}
a_{1,1}&...&a_{1,n+1}\\
...&...&...\\
a_{n+1,1}&...&a_{n+1,n+1}
\end{array}\right],$$
hence, due to \eqref{eq3.15}, follows expansion \eqref{eq3.12} for $n+1$. $\blacksquare$

Remark \ref{r2.3} and Theorems \ref{t3.1}, \ref{t3.2} imply that
\begin{equation}
r(\lambda)=\frac{(2i\lambda)^n}{\alpha_1...\alpha_n}(-1)^nb_1(\lambda^2)...b_n(\lambda^2)=r_1(\lambda)...r_n(\lambda)\label{eq3.16}
\end{equation}
where
\begin{equation}
r_k(\lambda)=-\frac{2i\lambda}{\alpha_k}b_k(\lambda^2)\quad(1\leq k\leq n).\label{eq3.17}
\end{equation}

\begin{theorem}\label{t3.3}
The scattering function $S(\lambda)$ has the multiplicative expansion
\begin{equation}
S(\lambda)=S_1(\lambda)...S_n(\lambda)\label{eq3.18}
\end{equation}
where
\begin{equation}
S_k(\lambda)=-\frac{r_k(-\lambda)}{r_k(\lambda)}\quad(1\leq k\leq n),\label{eq3.19}
\end{equation}
besides, $r_k(\lambda)$ are expressed via $b_k(\lambda^2)$ \eqref{eq3.13} by the formulas \eqref{eq3.17}.
\end{theorem}
\vspace{5mm}

{\bf 3.3} Multipliers $\{S_k(\lambda)\}$ \eqref{eq3.19} have natural interpretation.

{\bf The $n=1$ case.} The Jost solution $e_1(\lambda,x)$ equals
$$e_1(\lambda,x)=e^{i\lambda x}(\alpha_1\varphi_{2,1}(\lambda)-1)-\alpha_1\widetilde{v}_1^*(\lambda)\psi_1(\lambda,x),$$
due to \eqref{eq1.16}, \eqref{eq1.17}. Equations \eqref{eq1.19} and \eqref{eq1.13} imply
$$\alpha_1\varphi_{2,1}(\lambda)-1=\frac{\alpha_1}{2i\lambda}(\Phi_{1,1}^*(\lambda)-\Phi_{1,1}^*(-\lambda))-1=\frac{\alpha_1}{2i\lambda}[\Phi_{2,1}^*(\lambda)+\Phi_{1,1}(-\lambda)-
\widetilde{v}_1^*(\lambda)\widetilde{v}(\lambda)]$$
$$-1=-1-\frac{\alpha_1}{2i\lambda}[\widetilde{v}_1^*(\lambda)\widetilde{v}_1(-\lambda)-\Phi_{1,1}^*(\lambda)-\Phi_{1,1}(-\lambda)]+\frac{\alpha_1}{2i\lambda}\widetilde{v}_1^*(\lambda)[
\widetilde{v}_1(-\lambda)-\widetilde{v}(\lambda)]$$
$$=-(1+\alpha_1\langle R_0(\lambda^2)v_1,v_1\rangle)-\frac{\alpha_1}{2i\lambda}\widetilde{v}_1^*(\lambda)W_1(\lambda)$$
due to \eqref{eq3.7} and (see \eqref{eq3.4})
\begin{equation}
W_1(\lambda)\stackrel{\rm def}{=}\widetilde{v}_1(\lambda)-\widetilde{v}_1(-\lambda)=-2i\int\limits_{\mathbb{R}_+}\sin\lambda xv_1(x)dx,\label{eq3.20}\
\end{equation}
therefore
$$e_1(\lambda,x)=-e^{i\lambda x}(1+\alpha_1\langle R_2(\lambda^2)v_1,v_1\rangle)+\frac{\alpha_1}\lambda\widetilde{v}_1^*(\lambda)\left\{e^{i\lambda x}\int\limits_{\mathbb{R}_+}\sin\lambda tv_1(t)dt\right.$$
$$\left.-\int\limits_x^\infty\sin\lambda(t-\lambda)v_1(t)dt\right\}.$$
Using \eqref{eq3.3}, we arrive at the statement.

\begin{theorem}\label{t3.4}
The Jost solution $e_1(\lambda,x)$ \eqref{eq1.18} ($n=1$) is expressed via the boundary values on $\mathbb{R}$ from $\mathbb{C}_+$ of the resolvent $R_0(\lambda^2)$ \eqref{eq3.3} by the formula
\begin{equation}
e_1(\lambda,x)=-e^{i\lambda x}(1+\alpha_1\langle R_1(\lambda^2)v_1,v_1\rangle)+\alpha_1\widetilde{v}_1^*(\lambda)\cdot(R_1(\lambda^2)v_1)(x),\label{eq3.21}
\end{equation}
besides, $\langle e_1,v_1\rangle=-\langle e_0,v_1\rangle=-\widetilde{v}_1^*(\lambda)$ where $e_0(\lambda,x)=e^{i\lambda x}$ is the Jost solution of the operator $L_0$ \eqref{eq3.1}, \eqref{eq3.2}.
\end{theorem}

{\bf The $n\in\mathbb{N}$ case, ($n>1$).}

\begin{lemma}\label{l3.2}
Functions $\varphi_{k,s}(\lambda)$ \eqref{eq1.19} are expressed via the boundary values $\lambda+i0$ on $\mathbb{R}$ from $\mathbb{C}_+$ of the resolvent $R_0(\lambda^2)$ \eqref{eq3.3} by the formulas
\begin{equation}
\varphi_{k,s}(\lambda)=-\langle R_0(\lambda^2)v_k,v_s\rangle-\frac1{2i\lambda}\widetilde{v}_s^*(\lambda)W_k(\lambda)\quad(1\leq k,s\leq n)\label{eq3.22}
\end{equation}
where
\begin{equation}
W_k(\lambda)\stackrel{\rm def}{=}\widetilde{v}_k(\lambda)-\widetilde{v}_k(-\lambda)=-2i\int\limits_{\mathbb{R}_+}\sin\lambda xv_k(x)dx\quad(1\leq k\leq n).\label{eq3.23}
\end{equation}
\end{lemma}

P r o o f. Formulas \eqref{eq1.12}, \eqref{eq1.15} imply
$$\varphi_{k,s}(\lambda)=\frac1{2i\lambda}(\Phi_{k,s}^*(\lambda)-\Phi_{k,s}^*(-\lambda))=\frac1{2i\lambda}[\Phi_{k,s}^*(\lambda)+\Phi_{s,k}(-\lambda)-\widetilde{v}_s^*(\lambda)
\widetilde{v}_k(\lambda)]$$
$$=\frac1{2i\lambda}[\Phi_{k,s}^*(\lambda)+\Phi_{s,k}(-\lambda)-\widetilde{v}_s^*(\lambda)\widetilde{v}_k(-\lambda)]+\frac1{2i\lambda}\widetilde{v}_s^*(\lambda)[\widetilde{v}_k(-
\lambda)-\widetilde{v}_k(\lambda)],$$
which gives \eqref{eq3.22}, in view of \eqref{eq3.4}. $\blacksquare$

Using \eqref{eq1.16}, \eqref{eq1.17} and \eqref{eq3.22}, rewrite $e(\lambda,x)$ \eqref{eq1.18} as
$$e(\lambda,x)=(-1)^n\{e^{i\lambda x}$$
$$\times\det\left[
\begin{array}{cccc}
{\displaystyle
\begin{array}{ccc}
{\displaystyle\frac{\alpha_1}{2i\lambda}\widetilde{v}_1^*(\lambda)W_1(\lambda)}\\
+1+\alpha_1\langle R_0(\lambda^2)v_1,v_1\rangle
\end{array}}&...&{\displaystyle
\begin{array}{ccc}
{\displaystyle\frac{\alpha_n}{2i\lambda}\widetilde{v}_1^*(\lambda)W_n(\lambda)}\\
+\alpha_n\langle R_0(\lambda^2)v_n,v_1\rangle
\end{array}}\\
...&...&...\\
{\displaystyle
\begin{array}{cccccc}
\displaystyle{\frac{\alpha_1}{2i\lambda}\widetilde{v}_n^*(\lambda)W_1(\lambda)}\\
+\alpha_1\langle R_0(\lambda^2)v_1,v_n\rangle
\end{array}}&...&{\displaystyle
\begin{array}{ccc}
{\displaystyle\frac{\alpha_n}{2i\lambda}\widetilde{v_n^*}(\lambda)W_n(\lambda)}\\
+1+\alpha_n\langle R_0(\lambda^2)v_n,v_n\rangle
\end{array}}
\end{array}\right]$$
$$+\alpha_1\psi_1(\lambda,x)\cdot\det\left[
\begin{array}{ccccc}
\widetilde{v_1^*}(\lambda)&...&{\displaystyle\frac{\alpha_n}{2i\lambda}\widetilde{v}_1^*(\lambda)W_n(\lambda)+\alpha_n\langle R_0(\lambda^2)v_n,v_1\rangle}\\
...&...&...\\
\widetilde{v}_n^*(\lambda)&...&{\displaystyle\frac{\alpha_n}{2i\lambda}\widetilde{v}_n^*(\lambda)W_n(\lambda)+1+\alpha_n\langle R_0(\lambda^2)v_n,v_n\rangle}
\end{array}\right]+...$$
$$\left.+\alpha_n\psi_n(\lambda,x)\cdot\det\left[
\begin{array}{ccc}
{\displaystyle\frac{\alpha_1}{2i\lambda}\widetilde{v}_1^*(\lambda)W_1(\lambda)+1+\alpha_1\langle R_0(\lambda^2)v_1,v_1\rangle}&...&\widetilde{v}_1^*(\lambda)\\
...&...&...\\
{\displaystyle\frac{\alpha_1}{2i\lambda}\widetilde{v}_n^*(\lambda)W_1(\lambda)+\alpha_1\langle R_0(\lambda^2)v_1,v_n\rangle}&...&\widetilde{v}_n^*(\lambda)
\end{array}\right]\right\}.$$
Taking into account the proportionality of $\col[\widetilde{v}_1^*(\lambda)W_k(\lambda),...,\widetilde{v}_n^*(\lambda)W_k(\lambda)]$ ($1\leq k\leq n$), we obtain
$$e(\lambda,x)=(-1)^n\left\{e^{i\lambda x}\det\left[
\begin{array}{cccc}
\alpha_1\langle R_0(\lambda^2)v_1,v_1\rangle+1&...&\alpha_n\langle R_0(\lambda^2)v_n,v_1\rangle\\
...&...&...\\
\alpha_1\langle R_0(\lambda^2)v_1,v_n\rangle&...&\alpha_n\langle R_0(\lambda^2)v_n,v_n\rangle+1
\end{array}\right]\right.$$
$$+\alpha_1\left(\psi_1(\lambda,x)+\frac1{2i\lambda}e^{i\lambda x}W_1(\lambda)\right)\det\left[
\begin{array}{cccc}
\widetilde{v}_1^2(\lambda)&...&\alpha_n\langle R_0(\lambda^2)v_n,v_1\rangle\\
...&...&...\\
\widetilde{v}_n^*(\lambda)&...&\alpha_n\langle R_0(\lambda^2)v_n,v_n\rangle+1
\end{array}\right]+...$$
$$\left.+\alpha_n\left(\psi_n(\lambda,x)+\frac1{2i\lambda}e^{i\lambda x}W_n(\lambda)\right)\det\left[
\begin{array}{ccccc}
\alpha_1\langle R_0(\lambda^2)v_1,v_1\rangle+1&...&\widetilde{v}_1^*(\lambda)\\
...&...&...\\
\alpha_1\langle R_0(\lambda^2)v_1,v_n\rangle&...&\widetilde{v}_n^*(\lambda)
\end{array}\right]\right\}.$$

\begin{lemma}\label{l3.3}
For all $\lambda\in\mathbb{R}$, the following equalities are true:
$$\psi_k(\lambda,x)+\frac1{2i\lambda}e^{i\lambda x}W_k(\lambda)=-R_0(\lambda^2)v_k(x)\quad(1\leq k\leq n)$$
where $R_0(\lambda^2)$ are boundary values $\lambda+i0$ of resolvent \eqref{eq3.3} on $\mathbb{R}$ from $\mathbb{C}_+$.
\end{lemma}

Proof of the lemma is obvious,
$$\psi_k(\lambda,x)+\frac1{2i\lambda}e^{i\lambda x}W_k(\lambda)$$
$$=\frac1{2i\lambda}\left\{\int\limits_x^l(e^{i\lambda(t-x)}-e^{i\lambda(x-t)})v_k(t)+\int\limits_0^\infty(e^{i\lambda(x-t)}-e^{i\lambda(x+t)})v_k(t)dt\right\}$$
$$=-e^{i\lambda x}\int\limits_0^x\frac{\sin\lambda t}\lambda v_k(t)dt-\frac{\sin\lambda x}\lambda\int\limits_x^\infty e^{i\lambda t}v_k(t)dt=-R_0(\lambda^2)v_k(\lambda).$$
So,
\begin{equation}
e(\lambda,x)=(-1)^n\left\{e^{i\lambda x}\det\left[
\begin{array}{cccc}
\alpha_1\langle R_0(\lambda^2)v_1,v_1\rangle+1&...&\alpha_n\langle R_0(\lambda^2)v_n,v_1\rangle\\
...&...&...\\
\alpha_1\langle R_0(\lambda^2)v_1,v_n\rangle&...&\alpha_n\langle R_0(\lambda^2)v_n,v_n\rangle+1
\end{array}\right]\right.\label{eq3.24}
\end{equation}
$$-\alpha_1R_0(\lambda^2)v_1\cdot\det\left[
\begin{array}{ccc}
\widetilde{v}_1^*(\lambda)&...&\alpha_n\langle R_0(\lambda^2)v_n,v_1\rangle\\
...&...&...\\
\widetilde{v}_n^*(\lambda)&...&\langle\alpha_nR_0(\lambda^2)v_n,v_n\rangle+1
\end{array}\right]-...$$
$$\left.-\alpha_nR_0(\lambda^2)v_n\det\left[
\begin{array}{ccc}
\alpha_1\langle R_0(\lambda^2)v_1,v_1\rangle+1&...&\widetilde{v}_1^*(\lambda)\\
...&...&...\\
\alpha_1\langle R_0(\lambda^2)v_1,v_n\rangle&...&\widetilde{v}_n^*(\lambda)
\end{array}\right]\right\}.$$

\begin{theorem}\label{t3.5}
For the Jost solution $e_k(\lambda,x)$ of the operator $L_k$ \eqref{eq3.10} ($0\leq k\leq n$, $e_0(\lambda,x)=e^{i\lambda x}$ is the Jost solution of the operator $L_0$ \eqref{eq3.1}, \eqref{eq3.2}), the following recurrent formula is true:
\begin{equation}
e_k(\lambda,x)=-(1+\alpha_k\langle R_{k-1}(\lambda^2)v_k,v_k\rangle)e_{k-1}(\lambda,x)+\alpha_k\langle e_{k-1},v_k\rangle R_{k-1}(\lambda^2)v_k,\label{eq3.25}
\end{equation}
here $R_k(\lambda^2)=(L_k-\lambda^2I)^{-1}$ is resolvent of the operator $L_k$ \eqref{eq3.10} ($1\leq k\leq n$).
\end{theorem}

P r o o f. For $k=1$, the statement coincides with Theorem \ref{t3.4}. Using induction, we suppose that \eqref{eq3.25} holds for $k<n$ and prove that it is true for $k=n+1$. Prove that the function
$$\Phi(\lambda,x)=-(1+\alpha_n\langle R_{n-1}v_n,v_n\rangle)e_{n-1}(\lambda,x)+\alpha_n\langle e_{n-1},v_n\rangle R_{n-1}(\lambda^2)v_n$$
(the right-hand side of \eqref{eq3.25}) coincides, $e_n(\lambda,x)=e(\lambda,x)$ \eqref{eq3.24}. Using \eqref{eq3.17} and representation \eqref{eq3.25} for $e_{n-1}(\lambda,x)$, we obtain
$$\Phi(\lambda,x)=-\left\{1+\alpha_n\langle R_{n-1}(\lambda^2)v_n,v_n\rangle\right.$$
$$\left.-\frac{\alpha_{n-1}\langle R_{n-2}(\lambda^2)v_{n-1},v_n\rangle\alpha_n\langle R_{n-2}(\lambda^2)v_n,v_{n-1}\rangle}{1+\alpha_{n-1}\langle R_{n-2}(\lambda^2)v_{n-1},v_{n-1}\rangle}\right\}\cdot[-(1$$
$$+\alpha_{n-1}\langle R_{n-2}(\lambda^2)v_{n-1},v_{n-1}\rangle)e_{n-2}(\lambda,x)+\alpha_{n-1}\langle e_{n-2},v_{n-2}\rangle\langle R_{n-2}(\lambda^2)v_{n-1}]$$
$$+\alpha_n[-(1+\alpha_{n-1}\langle R_{n-2}(\lambda^2)v_{n-1},v_{n-1}\rangle\langle e_{n-2},v_n\rangle+\alpha_{n-1}\langle e_{n-2},v_{n-2}\rangle\langle R_{n-2}v_{n-1},v_n\rangle]$$
$$\times\left\{R_{n-2}(\lambda^2)v_n-\frac{\alpha_{n-1}\langle R_{n-2}(\lambda^2)v_n,v_{n-1}\rangle}{1+\alpha_{n-1}\langle R_{n-2}(\lambda^2)v_{n-1},v_{n-1}\rangle}R_{n-2}(\lambda^2)v_{n-1}\right\}.$$
Hence it follows that
$$\Phi(\lambda,x)=e_{n-2}(\lambda,x)\det\left[
\begin{array}{ccc}
\alpha_{n-1}\langle R_{n-2}(\lambda^2)v_{n-1},v_{n-1}\rangle+1&\alpha_n\langle R_{n-2}(\lambda^2)v_n,v_{n-1}\rangle\\
\alpha_{n-1}\langle R_{n-2}(\lambda^2)v_{n-1},v_n\rangle&\alpha_n\langle R_{n-2}(\lambda^2)v_n,v_n\rangle+1
\end{array}\right]$$
\begin{equation}
-\alpha_{n-1}R_{n-2}(\lambda^2)v_{n-1}\cdot\det\left[
\begin{array}{ccc}
\langle e_{n-2},v_{n-1}\rangle&\alpha_n\langle R_{n-2}(\lambda^2)v_n,v_{n-1}\rangle\\
\langle e_{n-2},v_n\rangle&\alpha_n\langle R_{n-2}(\lambda^2)v_n,v_n\rangle+1
\end{array}\right]\label{eq3.26}
\end{equation}
$$-\alpha_nR_{n-2}(\lambda^2)v_n\det\left[
\begin{array}{ccccc}
\alpha_{n-1}\langle R_{n-2}(\lambda^2)v_{n-1},v_{n-1}\rangle+1&\langle e_{n-2},v_{n-1}\rangle\\
\alpha_{n-1}\langle R_{n-2}(\lambda^2)v_{n-1},v_n\rangle&\langle e_{n-2},v_n\rangle
\end{array}\right].$$
Equality \eqref{eq3.26}, for $n=2$, coincides with \eqref{eq3.24}. Substituting representation \eqref{eq3.25} for $e_{n-2}(\lambda,x)$ into \eqref{eq3.26} and taking into account \eqref{eq3.14}, we arrive at equality \eqref{eq3.24} for $n=3$. Iterating this technique, we obtain that $\Phi(\lambda,x)=e(\lambda,x)$ \eqref{eq3.24}. $\blacksquare$

For the operator pair $\{L_1,L_0\}$, the Jost solution $e_1(\lambda,x)$ corresponding to the perturbed operator $L_1$ \eqref{eq3.10} is constructed \eqref{eq3.21} from the normalized Jost solution $e_2(\lambda,x)=e^{i\lambda x}$ ($e_0(\lambda,0)=1$) of the non-perturbed operator $L_0$ \eqref{eq3.1}, \eqref{eq3.2}. Analogously, for the pair $\{L_k,L_{k-1}\}$, we construct \eqref{eq3.25} the Jost solution $\widetilde{e}_k(\lambda,x)$ corresponding to the perturbed operator $L_k$ \eqref{eq3.10} by the normalized at zero Jost solution   $\widehat{e}_{k-1}(\lambda,x)$ ($=e_{k-1}(\lambda,x)/e_{k-1}(\lambda,0)$), $\widehat{e}_{k-1}(\lambda,0)=1$, of the non-perturbed operator $L_{k-1}$,
\begin{equation}
\widetilde{e}_{k-1}(\lambda,x)=-(1+\alpha_k\langle R_0(\lambda^2)v_k,v_k\rangle)\widehat{e}_{k-1}(\lambda,x)+\alpha_k\langle\widehat{e}_{k-1},v_k\rangle R_{k-1}(\lambda^2)v_k.\label{eq3.27}
\end{equation}

\begin{theorem}\label{t3.6}
Function $S_k(\lambda)$ \eqref{eq3.19} is the scattering coefficient of the pair $\{L_k,L_{k-1}\}$,
\begin{equation}
S_k(\lambda)=\widetilde{e}_k(-\lambda,0)/\widetilde{e}_k(\lambda,0)\label{eq3.28}
\end{equation}
where $\widetilde{e}_k(\lambda,x)$ is the Jost solution of the operator $L_k$ \eqref{eq3.10} calculated \eqref{eq3.27} from the normalized Jost solution $\widehat{e}_{k-1}(\lambda,x)$ ($\widehat{e}_{k-1}(\lambda,0)=1$) of the non-perturbed operator $L_{k-1}$.
\end{theorem}

Proof of the theorem follows from \eqref{eq3.17}, $\widehat{e}_{k-1}(\lambda,z)=-b_{k-1}(\lambda^2)$.

So, $S_k(\lambda)$ \eqref{eq3.19} is the scattering coefficient of the pair $\{L_k,L_{k-1}\}$ where $L_k=L_{k-1}+\alpha_k\langle.,v_k\rangle v_k$ is one-dimensional perturbation of the non-perturbed ({\bf background}) operator $L_{k-1}$, solution of which is normalized by identity at the point $x=0$.

\section{Inverse problem ($n=1$)}\label{s4}

This section is dedicated to the inverse problem for the pair $\{L_1,L_0\}$ where $L_0$ is an operator of the \eqref{eq3.1}, \eqref{eq3.2} kind, and $L_1=L_0+\alpha_1\langle.,v_1\rangle v_1$ \eqref{eq3.10}, and for $v_1$, \eqref{eq1.3} holds (see \cite{7}).

{\bf The case of $\alpha_1>0$}

{\bf 4.1} Equations \eqref{eq3.17}, \eqref{eq3.19} imply the boundary value Riemann problem \cite{17, 18},
\begin{equation}
S_1(\lambda)B_1^+(\lambda)=B_1^-(\lambda)\quad(\lambda\in\mathbb{R})\label{eq4.1}
\end{equation}
where $B_1^\pm(\lambda)$ are boundary values $\lambda\pm i0$ on $\mathbb{R}$ from $\mathbb{C}_\pm$ of the function $b_1(\lambda^2)=1+\alpha_1\langle R_0(\lambda^2)v_1,v_1\rangle$ \eqref{eq3.13}, besides,
\begin{equation}
B_1^\pm(\lambda)=1\pm\frac{\alpha_1}{2i\lambda}F_1(\pm\lambda)\quad(F_1(\lambda)=\widetilde{v}_1(-\lambda)\widetilde{v}_1^*(\lambda)-\Phi_{1,1}(\lambda)-\Phi_{1,1}(-\lambda)),
\label{eq4.2}
\end{equation}
due to \eqref{eq3.4}, \eqref{eq3.5}. Since (see \eqref{eq1.15}, \eqref{eq2.9})
$$\Re F_1(\lambda)=-\frac12|W_1(\lambda)|\quad(\forall\lambda\in\mathbb{R}, W_1(\lambda)=\widetilde{v}_1(\lambda)-\widetilde{v}_1(-\lambda))$$
and $F_1(\lambda)\in H_+^2$ (Remark \ref{r1.5}), then, using (iii) of Theorem \ref{t1.2}, we obtain
\begin{equation}
\begin{array}{ccccc}
{\displaystyle F_1(\lambda)=-\frac12|W_1(\lambda)|^2-\frac1{2\pi i}\int\limits_{\mathbb{R}}\hspace{-4.4mm}/\frac{dt}{t-\lambda}|W_1(t)|^2\left(=-\frac12|W_1(\lambda)|^2\right.}\\
{\displaystyle\left.-\frac\lambda{2\pi i}\int\limits_{\mathbb{R}_+}\hspace{-4.4mm}/\frac{dt}{t^2-\lambda^2}|W_1(t)|^2dt\right)\,(\lambda\in\mathbb{R}).}
\end{array}\label{eq4.3}
\end{equation}
and thus $F_1(\lambda)$ \eqref{eq4.3} is the boundary value $\lambda+i0$ of a Cauchy type integral \cite{17, 18},
\begin{equation}
{\mathcal F}_1(\lambda)=\frac1{2\pi i}\int\limits_{\mathbb{R}}\frac{dt}{t-\lambda}|W_1(t)|^2\quad(\lambda\in\mathbb{C}\setminus\mathbb{R})\label{eq4.4}
\end{equation}
($F_1(\lambda)=\mathcal{F}_1^+(\lambda)=\mathcal{F}_1(\lambda+i0)$). The boundary value $\mathcal{F}_1(\lambda-i0)=\mathcal{F}_1^-(\lambda)$ is
\begin{equation}
\mathcal{F}_1^-(\lambda)=\frac12|W_1(\lambda)|^2-\frac1{2\pi i}\int\limits_{\mathbb{R}}\hspace{-4.4mm}/\frac{dt}{t-\lambda}|W_1(t)|^2(=-F_1(-\lambda))\quad(\lambda\in\mathbb{R}).\label{eq4.5}
\end{equation}
So, \eqref{eq4.3}, \eqref{eq4.5} are Sokhoktsky formulas \cite{17, 18} for the Cauchy type integral \eqref{eq4.4}. Formulas \eqref{eq3.3}, \eqref{eq4.5} imply the expressions for $B_1^\pm(\lambda)$ \eqref{eq4.2}:
\begin{equation}
B_1^\pm(\lambda)=1\pm\frac{i\alpha_1}{4\lambda}|W_1(\lambda)|^2+\frac{\alpha_1}{4\pi\lambda}\int\limits_{\mathbb{R}}\hspace{-4.4mm}/\frac{dt}{t-\lambda}|W_1(t)|^2dt\quad(\lambda\in
\mathbb{R}),\label{eq4.6}
\end{equation}
besides,
\begin{equation}
B_1^\pm(\infty)=1.\label{eq4.7}
\end{equation}
Functions $B_1^\pm(\lambda)$ \eqref{eq4.6} coincide with the boundary values $\mathcal{B}_1(\lambda\pm i0)$ on $\mathbb{R}$ from $\mathbb{C}_\pm$ of the function
\begin{equation}
\mathcal{B}_1(\lambda)=1+\frac{\alpha_1}{2i\lambda}\mathcal{F}_1(\lambda)=1+\frac{\alpha_1}{4\lambda\pi}\int\limits_{\mathbb{R}}\frac{dt}{t-\lambda}|W_1(t)|^2\quad(\lambda\in\mathbb{C}
\setminus\mathbb{R})\label{eq4.8}
\end{equation}
where $\mathcal{F}(\lambda)$ is from \eqref{eq4.7}.

Representation $B_1^-(\lambda)=|B_1^-(\lambda)|e^{-i\zeta(\lambda)}$ ($\zeta(\lambda)=-\arg B_1^-(\lambda)$) and $B_1^+(\lambda)=\overline{B_1^-(\lambda)}$ ($\lambda\in\mathbb{R}$, see \eqref{eq4.6}) imply that
\begin{equation}
S_1(\lambda)=e^{-2i\zeta(\lambda)}\quad(\lambda\in\mathbb{R}).\label{eq4.9}
\end{equation}
Functions $B_1^\pm(\lambda)$ are holomorphic in $\mathbb{C}_\pm$ and don't vanish in $\mathbb{C}_\pm$ (Theorem \ref{t2.1}, $\alpha_1>0$). For $\lambda\in\mathbb{R}$, the functions $B_1^\pm(\lambda)$ can have $2q+1$ zeros (including $\lambda=0$, Lemma \ref{l2.7}), but since
$$B_1^\pm(0)=1+\frac{\alpha_1}{2\pi}\int\limits_{\mathbb{R}_+}\left|\frac{W_1(t)}t\right|^2dt>0\quad(\alpha_1>0),$$
then $\lambda=0$ is not a zero of $B_1^\pm(\lambda)$ and thus the set of real zeros of $B_1^\pm(\lambda)$ is $\{\pm\lambda_k\}_1^q$ ($0\leq q$).

If $q=0$, then $B^\pm(\lambda)$ don't have zeros in $\overline{\mathbb{C}^\pm}$ (index \cite{17, 18} of the Riemann problem \eqref{eq4.1} vanishes). Taking logarithm of \eqref{eq4.1}, $\ln B_1^+(\lambda)-\ln B_1^-(\lambda)=-\ln S_1(\lambda)$ ($\lambda\in\mathbb{R}$) and taking into account \eqref{eq4.9}, using Sokhotski formulas \cite{17, 18}, we obtain a solution to the boundary value problem,
\begin{equation}
\mathcal{B}_1(\lambda)=\exp\left\{\frac1\pi\int\limits_\mathbb{R}\frac{\zeta(t)}{t-\lambda}dt\right\}\quad(\lambda\in\mathbb{C}\setminus\mathbb{R}),\label{eq4.10}
\end{equation}
and it is unique due to \eqref{eq4.7}. Boundary values $B_1^\pm(\lambda)=\mathcal{B}_1(\lambda\pm i0)$ of the function $\mathcal{B}_1(\lambda)$ \eqref{eq4.10} are
\begin{equation}
\mathcal{B}_1^\pm(\lambda)=\exp\left\{\pm i\zeta(\lambda)+\frac1\pi\int\limits_\mathbb{R}\hspace{-4.4mm}/\frac{\zeta(t)}{t-\lambda}dt\right\}\quad(\lambda\in\mathbb{R})\label{eq4.11}
\end{equation}
and $B_1^-(\lambda)=\overline{B_1^+(\lambda)}$.

\begin{theorem}\label{t4.1}
If $\alpha_1>0$ ($n=1$) and $r(\lambda)$ doesn't have zeros in $\mathbb{R}\setminus\{0\}$, then the function $W_1(\lambda)=\widetilde{v}_1(\lambda)-\widetilde{v}_1(-\lambda)$ ($v_1(\lambda)$ is given by \eqref{eq1.8}) is expressed via $S_1(\lambda)$ \eqref{eq4.9} by the formula
\begin{equation}
|W_1(\lambda)|^2=\frac{4\lambda}{\alpha_1}\sin\zeta(\lambda)\cdot\exp\left\{\frac1\pi\int\limits_\mathbb{R}\hspace{-4.4mm}/\frac{\zeta(t)}{t-\lambda}dt\right\}(\lambda\in\mathbb{R}).
\label{eq4.12}
\end{equation}
\end{theorem}

Proof of the theorem follows from \eqref{eq4.5} and \eqref{eq4.11},
$$B_1^+(\lambda)-B_1^-(\lambda)=\frac{i\alpha_1}{2\lambda}|W_1(\lambda)|^2.$$
Normalization $\|v_1\|_{L^2}=1$ (see Remark \ref{r1.1}) and \eqref{eq4.12} unambiguously define $\alpha_1$. Since
\begin{equation}
W_1(\lambda)=-2i\int\limits_{\mathbb{R}_+}\sin\lambda xv_1(x)dx,\label{eq4.13}
\end{equation}
then the Parseval' equality \cite{13, 14} implies that ${\displaystyle\|W_1\|_{L^2}=4\frac2\pi\|v_1\|_{L^2}=\frac8\pi}$, therefore, upon integration of \eqref{eq4.12}, we obtain
\begin{equation}
\alpha_1=\frac\pi2\int\limits_{\mathbb{R}_+}\lambda\sin\zeta(\lambda)\exp\left\{\frac1\pi\int\limits_\mathbb{R}\hspace{-4.4mm}/\frac{\zeta(t)}{t-\lambda}dt\right\}d\lambda.\label{eq4.14}
\end{equation}

\begin{conclusion}\label{c1}
If conditions of Theorem \ref{t4.1} hold, then from $S_1(\lambda)=e^{-2i\zeta(\lambda)}$ \eqref{eq4.9}, we can retrieve: (a) number $\alpha_1$ \eqref{eq4.14}; (b) function $v_1(x)$ (ambiguously).
\end{conclusion}

So, from \eqref{eq4.12} by $S_1(\lambda)$ we calculate $W_1(\lambda)$ (ambiguously), and then apply the inverse Fourier sine transform to equality \eqref{eq4.13}.

\begin{remark}\label{r4.1}
Write \eqref{eq4.8} as
$$\frac{4\pi\lambda}{\alpha_1}(\mathcal{B}_1(\lambda)-1)=\int\limits_\mathbb{R}\frac{|W_1(t)|^2}{t-\lambda}dt\quad(\lambda\in\mathbb{C}\setminus\mathbb{R}$$
where the right-hand side is a Nevanlinna function. Applying the {\bf Perron -- Stieltjes formula} \cite{5, 12} to this function, we obtain
\begin{equation}
|W_1(x)|^2=\frac1\pi\lim\limits_{y\rightarrow+0}\frac{4\pi}{\alpha_1}\Im\{(x+iy)B_1^+(x+iy)\}=\frac{4\pi x}{\alpha_1}\Im B_1^+(x),\label{eq4.15}
\end{equation}
which gives \eqref{eq4.12} upon the substitution of $B_1^+(x)$ \eqref{eq4.11}.
\end{remark}
\vspace{5mm}

{\bf 4.2} Describe the class of functions $S_1(\lambda)$ that are the scattering coefficients of the pairs $\{L_1,L_0\}$ where $L_1=L_0+\alpha_1\langle.,v_1\rangle v_1$ and $v_1(x)$ satisfies condition \eqref{eq1.3}. It is more convenient to do this in terms of $\zeta(x)$. Relation $\zeta(\lambda)=\arg B_1^+(\lambda)$ ($\lambda\in\mathbb{R}$) and \eqref{eq4.6} imply
$$\tan\zeta(\lambda)=\frac{\alpha_1|W_1(\lambda)|^2}{\displaystyle4\lambda+\frac{\alpha_1}\pi\int\limits_{\mathbb{R}}\hspace{-4.4mm}/\frac{|W_1(t)|^2}{t-\lambda}dt}=\frac{\alpha_1|W_1
(\lambda)|^2}{\displaystyle2\lambda\left(2+\frac{\alpha_1}\pi\int\limits_{\mathbb{R}_+}\hspace{-4.4mm}/\frac{|W_1(t)|^2}{t^2-\lambda^2}dt\right)},$$
and thus
\begin{equation}
\zeta(\lambda)=\arctan\left\{\frac{\alpha_1|W_1(\lambda)|^2}{\displaystyle2\lambda\left(2+\frac{\alpha_1}\pi\int\limits_{\mathbb{R}_+}\hspace{-4.4mm}/\frac{|W_1(t)|^2}{t^2-\lambda^2}dt
\right)}\right\}\quad(\lambda\in\mathbb{R}).\label{eq4.16}
\end{equation}
Therefore $\zeta(\lambda)$ is a real continuous bounded odd ($\zeta(-\lambda)=-\zeta(\lambda)$) function, besides, $\zeta(\lambda)\geq0$ ($\lambda\in\mathbb{R}_+$) and $\zeta(\infty)=0$. Using $\arctan x\leq x$ ($\forall x\in\mathbb{R}_+$), from \eqref{eq4.16} we obtain
$$\zeta(\lambda)\leq\frac{\alpha_1|W_1(\lambda)|^2}{\displaystyle2\lambda\left(2+\frac{\alpha_1}\pi\int\limits_{\mathbb{R}_+}\hspace{-4.4mm}/\frac{|W_1(t)|^2}{t^2-\lambda^2)}dt\right)}
\quad(\forall\lambda\in\mathbb{R}_+),$$
therefore $\lambda\zeta(\lambda)\in L^1(\mathbb{R})$, which gives $\zeta(\lambda)\in L^2(\mathbb{R})$ and $\lambda\zeta(\lambda)\in L^2(\mathbb{R})$.

Function $(B_1^+(\lambda)-1)\lambda$ is bounded and continuous in $\overline{\mathbb{C}_+}$ (see \eqref{eq4.6} and Remark \ref{r1.5}), besides, $\lambda(B_1^+(\lambda)-1)\in L^2(\mathbb{R})$. Note that $\lambda(B^+(\lambda)-1)$ is differentiable and its derivative belongs to $L^2(\mathbb{R})$; $\{B_1^+(\lambda)-1+\lambda(B_1^+(\lambda))'\}\in L^2(\mathbb{R})$, and thus $\lambda(B_1^+(\lambda))'\in L^2(\mathbb{R})$ due to $(B_1^+(\lambda)-1)\in L^2(\mathbb{R})$. Hence and from \eqref{eq4.11} it follows that
$$\lambda(B_1^+(\lambda))'=\lambda\left(i\zeta'(\lambda)+\frac1\pi\int\limits_{\mathbb{R}}\hspace{-4.4mm}/\frac{\zeta(t)}{(t-\lambda)^2}dt\right)B_1^+(\lambda)\in L^2(\mathbb{R}),$$
and taking into account the boundedness of $(B_1^+(\lambda))^{-1}$, we obtain that
$$\lambda\left(i\zeta'(\lambda)+\frac1\pi\int\limits_{\mathbb{R}}\hspace{-4.4mm}/\frac{\zeta(t)}{(t-\lambda)^2}dt\right)=\lambda(B_1^+(\lambda))'(B_1^+(\lambda))^{-1}\in L^2(\mathbb{R}),$$
and thus
\begin{equation}
\lambda\zeta'(\lambda)\in L^2(\mathbb{R}),\quad\lambda\int\limits_\mathbb{R}\hspace{-4.4mm}/\frac{\zeta(t)}{(t-\lambda)^2}dt\in L^2(\mathbb{R}).\label{eq4.17}
\end{equation}
Differentiating \eqref{eq4.12}, we have
$$W_1(\lambda)(\overline{W_1(\lambda)})'+(|W_1(\lambda)|)'\overline{W_1(\lambda)}=\frac4{\alpha_1}\{\sin\zeta(\lambda)+\lambda\zeta'(\lambda)\cos\zeta(\lambda)$$
$$\left.+\lambda\sin\zeta(\lambda)\frac1\pi\int\limits_\mathbb{R}\hspace{-4.4mm}/\frac{\zeta(t)}{(t-\lambda)^2}dt\right\}\cdot\exp\left\{\frac1\pi\int\limits_\mathbb{R}\hspace{-4.4mm}/
\frac{\zeta(t)}{t-
\lambda}dt\right\},$$
and, taking into account \eqref{eq4.12},
\begin{equation}
(W_1(\lambda))'+\frac{W_1(\lambda)}{\overline{W_1(\lambda)}}(\overline{W_1(\lambda)})'=\left\{\frac1\lambda+\zeta'(\lambda)\cot\zeta(\lambda)+\frac1\pi\int\limits_{\mathbb{R}}
\hspace{-4.4mm}/
\frac{\zeta(t)dt}{(t-\lambda)^2}\right\}W_1(\lambda).\label{eq4.18}
\end{equation}

\begin{remark}\label{r4.2}
The set of real zeros $\lambda$ ($\in\mathbb{R}$) of the function $W_1(\lambda)$ ($W_1(\lambda)=0$) is at most countable.
\end{remark}

Therefore equality \eqref{eq4.18} holds for all $\lambda$ excluding, probably, a countable set. The left-hand side of \eqref{eq4.18} belongs to $L^2(\mathbb{R})$; in the right-hand side, the first and the third terms also belong to $L^2(\mathbb{R})$, consequently,
$$\zeta'(\lambda)\cot\zeta(\lambda)W_1(\lambda)\in L^2(\mathbb{R}),$$
hence, in view of \eqref{eq4.12},
\begin{equation}
\int\limits_\mathbb{R}(\zeta'(\lambda))^2\cot^2\zeta(\lambda)\lambda\sin\zeta(\lambda)d\lambda<\infty.\label{eq4.19}
\end{equation}

{\bf Class $\Omega_0$.} {\it A real continuous bounded odd function $\zeta(\lambda)$, such that $\zeta(\lambda)\geq0$ ($\forall\lambda\in\mathbb{R}_+$) is said to be of the {\bf class $\Omega_0$} if

(i) $\lambda\zeta(\lambda)\in L^1(\mathbb{R})$;

(ii) the function
$$F(\lambda)\stackrel{\rm def}{=}\int\limits_\mathbb{R}\hspace{-4.4mm}/\frac{\zeta(t)}{t-\lambda}dt\quad(\lambda\in\mathbb{R})$$
is bounded almost everywhere, differentiable and $\lambda F'(\lambda)\in L^2(\mathbb{R})$;

(iii) $\zeta(\lambda)$ has derivative almost everywhere and integral \eqref{eq4.19} converges.}

\begin{theorem}\label{t4.2}
Let $\zeta(\lambda)\in\Omega_0$ and, for $W_1(\lambda)$ \eqref{eq4.13}, equation \eqref{eq4.12} be true, then $v_1(x)$ satisfies condition \eqref{eq1.3}.
\end{theorem}
\vspace{5mm}

{\bf 4.3} Let $q\in\mathbb{N}$ and $\{\pm\lambda_k\}_1^q$ be the set of real zeros of $B_1^\pm(\lambda)$ \eqref{eq4.6}. Write the Riemann problem \eqref{eq4.1} as
\begin{equation}
\left(\frac{\lambda-i}{\lambda+i}\right)^{2q}S_1(\lambda)\frac{(\lambda+i)^{2q}}{\displaystyle\prod\limits_k(\lambda^2-\lambda_k^2)}B_1^+(\lambda)-\frac{(\lambda-i)^2q}{\displaystyle
\prod\limits_k(\lambda^2-\lambda_k^2)}B_1^-(\lambda)\quad(\lambda\in\mathbb{R}).\label{eq4.20}
\end{equation}
Functions
\begin{equation}
B_1^\pm(\lambda,i)\stackrel{\rm def}{=}\frac{(\lambda\pm i)^2q}{\displaystyle\prod\limits_k(\lambda^2-\lambda_k^2)}B_1^\pm(\lambda)\label{eq4.21}
\end{equation}
are holomorphic in $\mathbb{C}_\pm$ and don't have zeros in $\overline{\mathbb{C}}_\pm$, besides,
\begin{equation}
B_1^\pm(\infty,i)=1.\label{eq4.22}
\end{equation}
So, the Riemann problem \eqref{eq4.20} is
\begin{equation}
S_1(\lambda,i)B_1^+(\lambda,i)=B_1^-(\lambda,i)\quad(\lambda\in\mathbb{R})\label{eq4.23}
\end{equation}
where
\begin{equation}
S_1(\lambda,i)\stackrel{\rm def}{=}\left(\frac{\lambda-i}{\lambda+i}\right)^{2q}S_1(\lambda).\label{eq4.24}
\end{equation}
Unique (due to \eqref{eq4.22}) solution (see Subsection 4.1) to the boundary value problem \eqref{eq4.23} equals
\begin{equation}
\mathcal{B}_1(\lambda,i)=\exp\left\{-\frac1{2\pi i}\int\limits_\mathbb{R}\frac{\ln S_1(t,i)}{t-\lambda}dt\right\}\quad(\lambda\in\mathbb{C}\setminus\mathbb{R})\label{eq4.25}
\end{equation}
and its boundary values $\lambda\pm i0$ from $\mathbb{C}_\pm$ on $\mathbb{R}$ are
$$B_1^\pm(\lambda,i)=\exp\left\{\mp\ln S_1(\lambda,i)-\frac1{2\pi i}\int\limits_\mathbb{R}\frac{\ln S_1(t,i)}{t-\lambda}dt\right\}\quad(\lambda\in\mathbb{R}).$$
Using \eqref{eq4.21} and \eqref{eq4.24}, we have
\begin{equation}
B_1^\pm(\lambda)=\frac{\displaystyle\prod\limits_k(\lambda-\lambda_k^2)}{(\lambda^2+1)^q}\cdot\exp\left\{\mp\ln S_1(\lambda)-\frac1{2\pi i}\int\limits_\mathbb{R}\frac{\ln S_1(t,i)}{t-\lambda}dt\right\}\quad(\lambda\in\mathbb{R}).\label{eq4.26}
\end{equation}
Taking into account
$$\left|\frac{t-i}{t+i}\right|=1;\quad\arg\frac{t-i}{t+i}=-2\arccot\frac1t\quad(t\in\mathbb{R}),$$
we have
\begin{equation}
B_1^\pm(\lambda)=\psi^q(\lambda)\prod\limits_k(\lambda^2-\lambda_k^2)\exp\left\{\mp\ln S_1(\lambda)-\frac1{2\pi i}\int\limits_\mathbb{R}\frac{\ln S_1(t)}{t-\lambda}dt\right\}\quad(\lambda\in\mathbb{R})\label{eq4.27}
\end{equation}
where
\begin{equation}
\psi(\lambda)\stackrel{\rm def}{=}\frac1{\lambda^2+1}\exp\left\{\frac2\pi\int\limits_\mathbb{R}\hspace{-4.4mm}/\frac{\arccot1/t}{t-\lambda}dt\right\}.\label{eq4.28}
\end{equation}

\begin{theorem}\label{t4.3} Let $\alpha_1>0$ and $\{\pm\lambda_k\}_1^q$ ($\lambda_k>0$, $q\in\mathbb{N}$) be the set of the real zeros of the function $r(\lambda)$, then for $W_1(\lambda)$ the following equality holds:
\begin{equation}
|W_1(\lambda)|^2=\frac{4\lambda}{\alpha_1}\psi^q(\lambda)\prod\limits_k(\lambda^2-\lambda_k^2)\sin\zeta(\lambda)\cdot\exp\left\{\frac1\pi\int\limits_\mathbb{R}\hspace{-4.4mm}/\frac{
\zeta(t)}{t-\lambda}dt\right\}\quad(\lambda\in\mathbb{R})\label{eq4.29}
\end{equation}
where $S_1(\lambda)$ is given by \eqref{eq3.9} and $\psi(\lambda)$, correspondingly, \eqref{eq4.28}.
\end{theorem}

Number $\alpha_1$, due to normalization $\|v_1\|_{L^2}=1$, is explicitly calculated:
\begin{equation}
\alpha_1=\frac\pi2\int\limits_{\mathbb{R}_+}\lambda\psi^q(\lambda)\prod\limits_k(\lambda^2-\lambda_k^2)\sin\zeta(\lambda)\cdot\exp\left\{\frac1\pi\int\limits_{\mathbb{R}}\hspace{-4.4mm}/
\frac{\zeta(t)}{t-\lambda}dt\right\}d\lambda.\label{eq4.30}
\end{equation}

\begin{conclusion}\label{co2}
If suppositions of Theorem \ref{t4.3} hold, then from the data $\{S_1(\lambda),\{\pm\lambda_k\}_1^q\}$ we can recover: (a) number $\alpha_1$ \eqref{eq4.30}; (b) function $v_1(x)$ (ambiguously).
\end{conclusion}

Equations \eqref{eq4.25} and \eqref{eq4.26} imply that
\begin{equation}
B_1^+(\lambda)=Q_q(\lambda)\exp\left\{i\zeta(\lambda)+\frac1\pi\int\limits_\mathbb{R}\hspace{-4.4mm}/\frac{\zeta_q(t)}{t-\lambda}dt\right\}\quad(\lambda\in\mathbb{R})\label{eq4.31}
\end{equation}
where
\begin{equation}
Q_q(\lambda)\stackrel{\rm def}{=}\prod_1^q\frac{\lambda^2-\lambda_k^2}{\lambda^2+1};\quad\zeta_q(\lambda)\stackrel{\rm def}{=}\zeta(\lambda)+2q\arccot\frac1\lambda.\label{eq4.32}
\end{equation}
Analogously to Subsection 4.2, we obtain the description of scattering data.

{\bf Class $\Omega_q$.} {\it Let a set $\{\zeta(\lambda),\{\pm\lambda_k\}_1^q\}$ be given where $\lambda_k>0$ ($1\leq k\leq q$, $q\in\mathbb{N}$), $\zeta(\lambda)$ is a real bounded odd, continuous on $\mathbb{R}\setminus\{\pm\lambda_k\}_1^q$ function and $\zeta(\lambda)\geq0$ ($\forall\lambda\in\mathbb{R}_+$). This set is said to belong to the {\bf $\Omega_q$ class} if

(i) $\lambda\zeta(\lambda)\in L^1(\mathbb{R})$;

(ii) the function
$$F_q(\lambda)\stackrel{\rm def}{=}\int\limits_\mathbb{R}\hspace{-4.4mm}/\frac{\zeta_q(t)}{t-\lambda}dt$$
($\zeta_q(\lambda)$ is given by \eqref{eq4.32}) is bounded and has derivative almost everywhere, besides, $\lambda Q_q(\lambda)F'_q(\lambda)\in L^2(\mathbb{R})$ ($Q_q(\lambda)$ is from \eqref{eq4.32});

(iii) $\zeta(\lambda)$ is differentiable almost everywhere and
$$\int\limits_\mathbb{R}\lambda Q_q(\lambda)(\zeta'(\lambda))^2\cot^2\zeta(\lambda)\sin\zeta(\lambda)d\lambda<\infty.$$}

\begin{theorem}\label{t4.4}
If a set $\{\zeta(\lambda),\{\pm\lambda_k\}_1^q\}\in\Omega_q$ and $W_1(\lambda)$ is from \eqref{eq4.13}, then, if \eqref{eq4.29} is true, then $v_1(x)$ satisfies condition \eqref{eq1.3}.
\end{theorem}
\vspace{5mm}

{\bf The case of $\alpha_1<0$}

{\bf 4.4} Write $B_1^\pm(\lambda)$ \eqref{eq4.6} as
\begin{equation}
B_1^\pm(\lambda)=1\pm\frac{\alpha i}\lambda|W'_1(\lambda)|^2+\frac{\alpha_1}{2\pi}\int\limits_{\mathbb{R}_+}\hspace{-4.4mm}/\frac{|W_1(t)|^2}{t^2-\lambda^2}dt\quad(\lambda\in\mathbb{R}).\label{eq4.33}
\end{equation}
Function $B_1^+(\lambda)$ is holomorphically extendable into $\mathbb{C}_+$ and has pole $z_1=i\varkappa_1$ ($\varkappa_1>0$), therefore, in accordance with \eqref{eq2.27},
\begin{equation}
\int\limits_{\mathbb{R}_+}\frac{|W_1(t)|^2}{t^2+\varkappa_1^2}=-\frac{2\pi}{\alpha_1}>0.\label{eq4.34}
\end{equation}
Note that $\lambda=0$ is not a zero of the function $B_1^\pm(\lambda)$ \eqref{eq4.33}; really, using \eqref{eq4.31}, we obtain that
$$\left.\Re B_1^\pm(\lambda)\right|_{\lambda=0}=1+\frac{\alpha_1}{2\pi}\int\limits_{\mathbb{R}_+}\hspace{-4.4mm}/\frac{|W_1(t)|^2}{t^2}dt=\frac{\alpha_1\varkappa_1^2}{2\pi}\int\limits_{\mathbb{R}_+}
\hspace{-4.4mm}/\frac{|W_1(t)|^2}{t^2(t^2+\varkappa_1^2)}dt<0.$$
Using \eqref{eq4.34}, transform the appearance of function $\mathcal{B}_1(\lambda)$ \eqref{eq4.8},
\begin{equation}
\begin{array}{ccc}
{\displaystyle\mathcal{B}_1(\lambda)=\frac{\alpha_1}{2\pi}\left\{\int\limits_{\mathbb{R}_+}\frac{|W_1(t)|^2}{t^2-\lambda^2}dt-\int\limits_{\mathbb{R}_+}\frac{|W_1(t)|^2}{t^2+\varkappa_1^2
}dt\right\}=\frac{\alpha_1}{2\pi}(\lambda^2+\varkappa_1^2)\int\limits_{\mathbb{R}_+}\frac{|W_1(t)|^2}{t^2-\lambda^2}\frac{dt}{t^2+\varkappa_1^2}}\\
{\displaystyle=\frac{\alpha_1}{4\pi\lambda}(\lambda^2+\varkappa_1^2)\int\limits_\mathbb{R}\frac{|W_1(t)|^2}{t-\lambda}\frac{dt}{t^2+\varkappa_1^2}\quad(\lambda\in\mathbb{C}\setminus
\mathbb{R}).}
\end{array}\label{eq4.35}
\end{equation}
The boundary values $B_1^\pm(\lambda)=\mathcal{B}_1(\lambda\pm i0)$ on $\mathbb{R}$ from $\mathbb{C}_\pm$ of the function $\mathcal{B}_1(\lambda)$ \eqref{eq4.35} are
\begin{equation}
B_1^\pm(\lambda)=\pm\frac{\alpha_1i}{4\lambda}|W_1(\lambda)|^2+\frac{\alpha_1}{4\lambda\pi}(\lambda^2+\varkappa_1)^2\int\limits_\mathbb{R}\hspace{-4.4mm}/\frac{|W_1(t)|^2}{t-\lambda}
\frac{dt}{t^2+\varkappa_1^2}\quad(\lambda\in\mathbb{R}),\label{eq4.36}
\end{equation}
and thus
\begin{equation}
B_1^+(\lambda)-B_1^-(-\lambda)=\frac{\alpha i}{2\lambda}|W_1(\lambda)|^2.\label{eq4.37}
\end{equation}
Functions $B_1^\pm(\lambda)$ vanish when $\lambda=\pm z_1=\pm i\varkappa$ ($\varkappa>0$) and on the real axis at the points $\{\pm\lambda_k\}_1^q$ ($\lambda_k>0$). Rewrite the boundary value problem \eqref{eq4.1} as
\begin{equation}
S_1(\lambda,i,\varkappa_1)B_1^+(\lambda,i,\varkappa)=B_1^-(\lambda,i,\varkappa)\label{eq4.38}
\end{equation}
where
\begin{equation}
B_1^\pm(\lambda,i,\varkappa)=\frac{(\lambda\pm i)^{iq+1}}{\displaystyle\prod\limits_k(\lambda^2-\lambda_k^2)}\frac{B_1^\pm(\lambda)}{(\lambda\mp i\varkappa)};\quad S_1(\lambda,i,\varkappa)=S_1(\lambda)\left(\frac{\lambda-i}{\lambda+i}\right)^{iq+1}\frac{\lambda-i\varkappa}{\lambda+i\varkappa}.\label{eq4.39}
\end{equation}
Functions $B_1^\pm(\lambda,i,\varkappa)$ are holomorphically extendable into $\mathbb{C}_\pm$ and don't have zeros in $\overline{\mathbb{C}_\pm}$, besides,
\begin{equation}
B_1^\pm(\infty,i,\varkappa)=1.\label{eq4.40}
\end{equation}
Index of the Riemann problem \eqref{eq4.38} vanishes, and its unique solution (see \eqref{eq4.40}) equals
$$\mathcal{B}_1(\lambda,i,\varkappa)=\exp\left\{-\frac1{2\pi i}\int\limits_\mathbb{R}\frac{\ln S_1(t,i,\varkappa)}{t-\lambda}\right\}\quad(\lambda\in\mathbb{C}\setminus\mathbb{R}),$$
besides,
$$B_1^\pm(\lambda,i,\varkappa)=\mathcal{B}_1(\lambda\pm i0,i,\varkappa)$$
$$=\exp\left\{\mp\frac12\ln S_1(\lambda,i,\varkappa)-\frac1{2\pi i}\int\limits_\mathbb{R}\hspace{-4.4mm}/\frac{\ln S_1(t,i,\varkappa)}{t-\lambda}dt\right\}\,(\lambda\in\mathbb{R}),$$
and, in view of \eqref{eq4.40}, hence we find
$$B_1^\pm(\lambda)=\sqrt{\lambda^2+\varkappa^2}\frac{\displaystyle\prod\limits_k(\lambda^2-\lambda_k^2)}{(\lambda^2+1)^{q+1/2}}\exp\left\{\mp\frac12\ln S_1(\lambda)-\frac1{2\pi i}\int\limits_\mathbb{R}\hspace{-4.4mm}/\frac{\ln S_1(t,i,\varkappa)}{t-\lambda}dt\right\}$$
$$=(\lambda^2+\varkappa^2)\prod\limits_k(\lambda^2-\lambda_k^2)\psi^{q+1/2}(\lambda)\psi^{1/2}(\lambda,\varkappa)\exp\left\{\mp\frac12\ln S_1(\lambda)-\frac1{2\pi i}\int\limits_\mathbb{R}\hspace{-4.4mm}/\frac{\ln S_1(t)}{t-\lambda}dt\right\}$$
$$(\lambda\in\mathbb{R})$$
where
$$\psi(\lambda,\varkappa)\stackrel{\rm def}{=}\frac1{\lambda^2+\varkappa^2}\exp\left\{\frac1\pi\int\limits_{\mathbb{R}}\frac{\arccot\varkappa/t}{t-\lambda}dt\right\}.$$
It is easy to show that $\psi(\lambda,\varkappa)$ does not depend on $\varkappa$, therefore $\psi(\lambda,\varkappa)=\psi(\lambda,1)=\psi(\lambda)$, and thus
\begin{equation}
B_1^\pm(\lambda)=(\lambda^2+\varkappa^2)\prod\limits_k(\lambda^2-\lambda_k^2)\psi^{q+1}(\lambda)\cdot\exp\left\{\pm i\zeta(\lambda)+\frac1\pi\int\limits_{\mathbb{R}}\hspace{-4.4mm}/\frac{\zeta(t)}{t-\lambda}dt\right\}\,(\lambda\in\mathbb{R}).\label{eq4.41}
\end{equation}
Using \eqref{eq4.37}, we arrive at the theorem.

\begin{theorem}\label{t4.5}
Let $\alpha_1<0$ and $z=i\varkappa_1$ ($\varkappa_1>0$) be a zero of $r(\lambda)$ from $\mathbb{C}_+$ and $\{0,,\pm\lambda_k\}_1^q$ ($\lambda_k>0$, $q\in\mathbb{N}$) be the set of real zeros of $r(\lambda)$, then for $W_1(\lambda)$ the following representation is true:
\begin{equation}
|W_1(\lambda)|^2=\frac{4\lambda}{\alpha_1}\psi^{q+1}(\lambda)(\lambda^2+\varkappa^2)\prod\limits_k(\lambda^2-\lambda_k^2)\sin\zeta(\lambda)\cdot\exp\left\{\frac1\pi\int\limits_\mathbb{R}
\hspace{-4.4mm}/\frac{\zeta(t)}{t-\lambda}dt\right\}\,(\lambda\in\mathbb{R})\label{eq4.42}
\end{equation}
where $S_1(\lambda)$ is given by \eqref{eq4.9} and $\psi(\lambda)$, correspondingly, by \eqref{eq4.22}.
\end{theorem}

Number $\alpha_1$ is calculated from \eqref{eq4.42} ($\|v_1\|_{L^2}=1$),
\begin{equation}
\alpha_1=\frac\pi2\int\limits_{\mathbb{R}_+}\psi^{q+1}(\lambda)(\lambda^2+\varkappa^2)\prod\limits_k(\lambda^2-\lambda_k^2)\sin\zeta(\lambda)\cdot\exp\left\{\frac1\pi
\int\limits_\mathbb{R}\hspace{-4.4mm}/\frac{\zeta(t)}{t-\lambda}dt\right\}d\lambda.\label{eq4.43}
\end{equation}

\begin{conclusion}\label{co3}
If conditions of Theorem \ref{t4.5} hold, then from the data $\{S_1(\lambda),\{\pm\lambda_k\},$ $\varkappa\}$ we can restore (a) number $\alpha_1$ \eqref{eq4.43}; (b) function $v_1(x)$ (ambiguously).
\end{conclusion}

Description of the scattering data when $\alpha_1<0$ lies in the following.

{\bf Class $\Omega_{q,\varkappa}$.} {\it Let there be given a totality $\{\zeta(\lambda),\{\pm\lambda_k\}_1^q,\varkappa\}$ where $\lambda_k>0$ ($1\leq k\leq q$, $q\in\mathbb{N}$); $\varkappa>0$; $\zeta(\lambda)$ is a real, bounded, odd, continuous on $\mathbb{R}\setminus\{\pm\lambda_k\}_1^q$ function and $\zeta(\lambda)>0$ ($\forall\lambda\in\mathbb{R}_+$). This set is said to belong to class $\Omega_{q,\varkappa}$ if

(i) $\lambda\zeta(\lambda)\in L^1(\mathbb{R})$;

(ii) function
$$F_{q+1}(\lambda)\stackrel{\rm def}{=}\int\limits_\mathbb{R}\hspace{-4.4mm}/\frac{\zeta_{q+1}(t)}{t-\lambda}dt$$
is bounded and almost everywhere differentiable, besides, $\lambda Q_{q,\varkappa}(\lambda)F'_q(\lambda)\in L^2(\mathbb{R})$ where
\begin{equation}
Q_{q,\varkappa}(\lambda)\stackrel{\rm def}{=}\frac{\lambda^2+\varkappa^2}{\lambda^2+1}\prod\limits_k\frac{\lambda^2-\lambda_k^2}{\lambda^2+1};\quad\zeta_{q+1}(\lambda)\stackrel{\rm def}{=}\zeta(\lambda)+2(q+1)\arccot1/\lambda;\label{eq4.44}
\end{equation}

(iii) function $\zeta(\lambda)$ has derivative almost everywhere and
$$\int\limits_\mathbb{R}\lambda Q_{q,\varkappa}(\lambda)(\zeta'(\lambda))^2\cot^2\zeta(\lambda)\sin\zeta(\lambda)d\lambda<\infty.$$}

\begin{theorem}\label{t4.6}
If $\{\zeta(\lambda),\{\pm\lambda_k\}_1^q,\varkappa\}\in\Omega_{q,\varkappa}$, and $W_1(\lambda)$ is given by \eqref{eq4.13} and equality \eqref{eq4.42} holds, then $v_1(x)$ satisfies condition \eqref{eq1.3}.
\end{theorem}

The inclusions $\Omega_0\subseteq\Omega_q\subseteq\Omega_{q,\varkappa}$ take place.

\begin{conclusion}\label{co4}
Inverse scattering problem for the pair $\{L_1,L_0\}$ has solution. From the scattering data $\{\zeta(\lambda),\{\pm\lambda_k\}_1^q\}\in\Omega_q$ ($\alpha_1>0$) or $\{\zeta(\lambda),\{\pm\lambda_k\}_1^q,\varkappa\}\in\Omega_{q,\varkappa}$ ($\alpha_1<0$), we can recover number $\alpha_1$ unambiguously and function $v_1(x)$ ambiguously.
\end{conclusion}

\section{Inverse problem ($n\geq2$)}\label{s5}

{\bf The case of $n=2$}

{\bf 5.1} Scattering function $S(\lambda)$ of the pair $\{L_2,L_0\}$ (${\displaystyle L_2=L_0+\sum\limits_1^2\alpha_k\langle.,v_k\rangle v_k}$), due to \eqref{eq3.18}, equals
\begin{equation}
S(\lambda)=S_1(\lambda)S_2(\lambda)\label{eq5.1}
\end{equation}
where $\{S_k(\lambda)\}$ are given by \eqref{eq3.19}.

Multiplier $S_1(\lambda)$ is the coefficient of the Riemann boundary value problem \eqref{eq4.1},
\begin{equation}
S_1(\lambda)B_1^+(\lambda)=B_1^-(\lambda)\quad(\lambda\in\mathbb{R})\label{eq5.2}
\end{equation}
where $B_1^\pm(\lambda)$ are boundary values on $\mathbb{R}$ from $\mathbb{C}_\pm$ of the function $b_1(\lambda^2)=1+\alpha_1\langle R_0(\lambda^2)v_1,v_1\rangle$ \eqref{eq3.13} (or $\mathcal{B}_1(\lambda)$ \eqref{eq4.8}). From the scattering data $\{\zeta_1(\lambda),\{\pm\lambda_k(1)\}_1^{q_1},$ $\varkappa_1\}\in\Omega_{q,\varkappa}$ ($S_1(\lambda)=\exp(-2i\zeta_1(\lambda))$ \eqref{eq4.9}), we can unambiguously calculate (see Section \ref{s4}) number $\alpha_1\in\mathbb{R}$ and function $|W_1(\lambda)|^2$ ($W_1(\lambda)=\widetilde{v}_1(\lambda)-\widetilde{v}_1(-\lambda)$; $\lambda\in\mathbb{R}$), whence $v_1(x)$ is (ambiguously) defined. So, from the set $\{\zeta_1(\lambda),\{\pm\lambda_k(1)\}_1^{q_1},\varkappa_1\}$ we can (ambiguously) find the operator
\begin{equation}
L_1=L_0+\alpha_1\langle.,v_1\rangle v_1.\label{eq5.3}
\end{equation}
Function $S_2(\lambda)$ is the coefficient of the boundary value problem
\begin{equation}
S_2(\lambda)B_2^+(\lambda)=B_2^-(\lambda)\quad(\lambda\in\mathbb{R})\label{eq5.4}
\end{equation}
where $B_2^\pm(\lambda)$ are boundary $\lambda+i0$ values of the function $b_2(\lambda^2)=1+\alpha_2\langle R_1(\lambda^2)v_2,v_2\rangle$ \eqref{eq3.13} and $R_1(z)=(L_1-zI)^{-1}$ is the resolvent of operator $L_1$ \eqref{eq5.3}. Number $\alpha_2$ is calculated from $b_2(\lambda^2)$ due to normalization $\|v_2\|=1$. Equation \eqref{eq3.12} yields
\begin{equation}
b_1(z)b_2(z)=\det(I+\alpha T(z))=\det\left[
\begin{array}{ccc}
1+\alpha_1\langle R_0(z)v_1,v_1\rangle&\alpha_1\langle R_0(z)v_2,v_1\rangle\\
\alpha_2\langle R_0(z)v_1,v_2\rangle&1+\alpha_2\langle R_0(z)v_2,v_2\rangle
\end{array}\right]\label{eq5.5}
\end{equation}
where $\alpha$ and $T(z)$ are given by \eqref{eq1.25} and \eqref{eq3.6} ($n=2$). As a result, we arrive at the following problem of extension of a scalar Nevanlinna function up to a $2\times2$ matrix-valued Nevanlinna function.

{\bf Problem of $N$-extension.} {\it Consider a self-adjoint operator $L$ in a Hilbert space $H$, and let
\begin{equation}
A_1(z)\stackrel{\rm def}{=}\beta_1+\langle R(z)v_1,v_1\rangle;\quad A_2(z)\stackrel{\rm def}{=}\det\left[
\begin{array}{ccc}
\beta_1+\langle R(z)v_1,v_1\rangle&\langle R(z)v_1,v_2\rangle\\
\langle R(z)v_2,v_1\rangle&\beta_2+\langle R(z)v_2,v_2\rangle
\end{array}\right]\label{eq5.6}
\end{equation}
where $\beta_k\in\mathbb{R}$ ($\beta_k\not=0$, $k=1$, $2$); $\{v_k\}_1^2$ are linearly independent vectors from $H$; $R(z)=(L-zI)^{-1}$. How to define the function $\beta_2+\langle R(z)v_2,v_2\rangle$ from the set $\{A_k(z)\}_1^2$ and $\{\beta_k\}_1^2$? Describe the technique of finding this function and characterize the degree of ambiguity of the problem solution.}

Equation \eqref{eq3.7} implies that
\begin{equation}
\mathcal{B}(\lambda)\stackrel{\rm def}{=}I+\alpha T(\lambda^2)=I-\frac\alpha{2i\lambda}F(\lambda)\quad(\lambda\in\mathbb{C}\setminus\mathbb{R})\label{eq5.7}
\end{equation}
where $F(\lambda)$ is from \eqref{eq2.8}. Given \eqref{eq2.10}, we find the boundary $\lambda\pm i0$ values on $\mathbb{R}$ from $\mathbb{C}_\pm$ of the function $F(\lambda)$,
$$F^\pm(\lambda)=\pm\frac12W^+(\lambda)W(\lambda)+\frac1{2\pi i}\int\limits_\mathbb{R}\hspace{-4.4mm}/\frac{dt}{t-\lambda}W^+(t)W(t)\quad(\lambda=\lambda\pm i0\in\mathbb{R}),$$
and hence $F^\pm(\lambda)=\mathcal{F}(\lambda\pm i0)$ are boundary values of the Cauchy type integral \cite{7, 15},
$$\mathcal{F}(\lambda)\stackrel{\rm def}{=}\frac1{2\pi i}\int\limits_\mathbb{R}\frac{dt}{t-\lambda}W^+(t)W(t)\quad(\lambda\in\mathbb{C}\setminus\mathbb{R})$$
with the $(2\times2)$ matrix-valued density $W^+(t)W(t)$ ($\rank W^+(t)W(t)=1$). Therefore, in view of \eqref{eq5.7},
\begin{equation}
\mathcal{B}(\lambda)=I+\frac\alpha{4\pi\lambda}\int\limits_\mathbb{R}\frac{dt}{t-\lambda}W^+(t)W(t)\quad(\lambda\in\mathbb{C}\setminus\mathbb{R}),\label{eq5.8}
\end{equation}
besides,
\begin{equation}
\mathcal{B}(\lambda\pm i0)=I+\frac\alpha{4\lambda}\left\{\pm iW^+(\lambda)W(\lambda)+\frac1\pi\int\limits_\mathbb{R}\hspace{-4.4mm}/\frac{dt}{t-\lambda}W^+(t)W(t)\right\}\quad(\lambda\in\mathbb{R}).\label{eq5.9}
\end{equation}
Since
\begin{equation}
W^+(t)W(t)=\left[
\begin{array}{ccc}
W_1^+(t)W_1(t)&W_1^+(t)W_2(t)\\
W_2^+(t)W_1(t)&W_2^+(t)W_2(t)
\end{array}\right],\label{eq5.10}
\end{equation}
here
\begin{equation}
W_k(\lambda)=\widetilde{v}_k(\lambda)-\widetilde{v}_k(-\lambda)=-2i\int\limits_{\mathbb{R}_+}\sin\lambda xv_k(x)dx\quad(1\leq k\leq2),\label{eq5.11}
\end{equation}
then function $A_2(\lambda)$ equals
\begin{equation}
A_2(\lambda)=\det\left[
\begin{array}{ccc}
{\displaystyle\beta_1+\frac1{4\pi\lambda}\int\limits_\mathbb{R}\frac{dt}{t-\lambda}|W_1(\lambda)|^2}&{\displaystyle\frac1{4\pi\lambda}\int\limits_\mathbb{R}\frac{dt}{t-\lambda}
\overline{W}_1(t)W_2(t)}\\
{\displaystyle\frac1{4\pi\lambda}\int\limits_\mathbb{R}\frac{dt}{t-\lambda}\overline{W_2(t)}W_1(t)}&{\displaystyle\beta_2+\frac1{4\pi\lambda}\int\limits_\mathbb{R}\frac{dt}{t-\lambda}|
W_2(t)|^2}
\end{array}\right]\,(\lambda\in\mathbb{C}\setminus\mathbb{R}).\label{eq5.12}
\end{equation}
As a result, problem of $N$ expansion is reduced to the finding the function $|W_2(t)|^2$ from the known data $\{\beta_k\}$, $|W_1(t)|^2$, $A_2(\lambda)$ \eqref{eq5.12}.

\begin{remark}\label{r5.1}
Knowing the solution to the problem of $N$-expansion, similar to considerations of Section \ref{s4}, we, from the function $1+\alpha_2\langle R_0(\lambda^2)v_2,v_2\rangle$, construct the scattering data $\{\zeta_2(\lambda),\{\pm\lambda_k(2)\}_1^{q_2},\varkappa_2\}\in\Omega_{q,\varkappa}$ ($\exp(-2i\zeta_2(\lambda))$ is the coefficient of the boundary value problem for the function $1+\alpha_2\langle R_0(\lambda^2)v_2,v_2\rangle$), using which we reconstruct (Section \ref{s4}) number $\alpha_2$ and function $v_2(x)$ (ambiguously).
\end{remark}
\vspace{5mm}

{\bf 5.2} We will need the following definitions.

\begin{definition}\label{d5.1}
Let $L$ be a self-adjoint operator in a Hilbert space $H$, and $E_k$ be its resolution of identity. By $G(h)$, we denote the subspace corresponding to a vector $h\in H$,
\begin{equation}
G(h)\stackrel{\rm def}{=}\span\{E_th:t\in\mathbb{R}\}.\label{eq5.13}
\end{equation}
A vector $f\in H$ is said to be {\bf $L$-orthogonal to the vector $h$} if $f\perp G(h)$.
\end{definition}

\begin{definition}\label{d5.2}
The kernel $K_n(x,t)=\sum\limits_{k=1}^n\alpha_kv_k(x)\overline{v_k(t)}$ \eqref{eq1.1} is said to be of the {\bf class of $L_0$-orthogonal kernels} ($L_0$ is given by \eqref{eq3.1}) if the function $v_k(x)$ is $L_0$-orthogonal to the functions $v_s(x)$ ($1\leq s\leq k-1$) for all $k$ ($1\leq k\leq n$).
\end{definition}

Let $n=2$ and the kernel $K_2(x,t)$ be $L_0$-orthogonal. Then $v_2\perp G(v_1)$ and thus $\left.(L_1-L_0)\right|_{G^\perp(v_1)}=0$, therefore $\langle R_0(z)v_1,v_2\rangle=0$ and matrix $T(z)$ is diagonal. In this case,
\begin{equation}
b_2(\lambda^2)=1+\alpha_2\langle R_0(\lambda^2)v_2,v_2\rangle,\label{eq5.14}
\end{equation}
and $S_2(\lambda)$ is the coefficient of the Riemann problem
\begin{equation}
S_2(\lambda)B_2^+(\lambda)=B_2^+(\lambda)\quad(\lambda\in\mathbb{R})\label{eq5.15}
\end{equation}
where $B_2^\pm(\lambda)$ are boundary $\lambda+i0$ values on $\mathbb{R}$ from $\mathbb{C}_\pm$ of the function $b_2(\lambda^2)$ \eqref{eq5.14}. Following Section \ref{s4}, using the data $\{\zeta_2(\lambda),\{\pm\lambda_k(2)\}_1^{q_2},\varkappa_2\}\in\Omega_{q,\varkappa}$ ($s_2(\lambda)=\exp(-2i\zeta_2(\lambda))$), the inverse problem is solved and number
$\alpha_2\in\mathbb{R}$ is found, and also, function $v_2(x)$ is ambiguously found.

\begin{theorem}\label{t5.1}
If $n=2$ and the kernel $K_2(\lambda,t)$ \eqref{eq1.1} belongs to the class of $L_0$-orthogonal kernels ($L_0$ is given by \eqref{eq3.1}), then inverse problem is solved using the Riemann problems \eqref{eq5.2} and \eqref{eq5.15} where $\{S_k(\lambda)\}$ are multipliers from the multiplicative expansion \eqref{eq5.1} of the scattering function $S(\lambda)$ of the pair $\{L_2,L_0\}$.
\end{theorem}
\vspace{5mm}

{\bf The case of $n\in\mathbb{N}$}

{\bf 5.3} Equation \eqref{eq3.18} implies that the scattering function $S(\lambda)$ of the pair $\{L_n,L_0\}$ (${\displaystyle L_n=L_0+\sum\limits_1^n\alpha_k\langle.,v_k\rangle v_k}$ \eqref{eq3.10}) equals
\begin{equation}
S(\lambda)=S_1(\lambda)S_2(\lambda)...S_n(\lambda)\label{eq5.16}
\end{equation}
where $\{S_k(\lambda)\}$ are from \eqref{eq3.19}. If the kernel $K_n(x,t)$ \eqref{eq1.1} belongs to the class of $L_0$-orthogonal kernels, then for each of the multipliers $S_k(\lambda)$ in \eqref{eq5.16} we have the Riemann boundary value problem
\begin{equation}
S_k(\lambda)B_k^+(\lambda)=B_k^+(\lambda)\quad(1\leq k\leq n,\lambda\in\mathbb{R})\label{eq5.17}
\end{equation}
where $B_k^\pm(\lambda)$ are boundary $\lambda\pm i0$ values on $\mathbb{R}$ from $\mathbb{C}_\pm$ of the function
\begin{equation}
b_k(\lambda^2)=1+\alpha_k\langle R_0(\lambda^2)v_k,v_k\rangle\quad(1\leq k\leq n).\label{eq5.18}
\end{equation}
Following Section \ref{s4}, from the data $\{\zeta_k(\lambda),\{\pm\lambda_p(k)\}_1^{q_k},\varkappa_k\}\in\Omega_{q,\varkappa}$ ($S_k(\lambda)=\exp\{-2i\zeta_k(\lambda)\}$), we can calculate the real number $\alpha_k$ and a function $v_k(x)$ (ambiguously).

\begin{theorem}\label{t5.2}
If the kernel $K_n(x,t)$ \eqref{eq1.1} of the operator $L_n$ \eqref{eq3.10} belongs to the class of $L_0$-orthogonal kernels ($L_0$ is given by \eqref{eq3.1}), then the inverse problem is solvable and the numbers $\alpha_k$ and $|W_k(\lambda)|^2$ ($W_k(\lambda)$ are from \eqref{eq5.11}) are defined by the scattering data $\{\zeta_k(\lambda),\{\pm\lambda_p(k)\}_1^{q_k},\varkappa_k\}\in\Omega_{q,\varkappa}$ ($S_k(\lambda)=\exp(-2i\zeta_k(\lambda)$) via the formulas \eqref{eq4.42} and \eqref{eq4.43}.
\end{theorem}

So, an inverse problem is solvable in the case when the scattering coefficient $S(\lambda)$ has the multiplicative expansion \eqref{eq5.16}, besides, each of multipliers $S_k(\lambda)$ belongs to the class $\Omega_{q_k,\varkappa_k}$ (see Section \ref{s4}).

\renewcommand{\refname}{ \rm 
\end{Large}
\end{document}